\documentclass[11pt]{article}
\usepackage{latexsym}
\usepackage[all]{xy}
\usepackage{amsmath}
\usepackage{amssymb}
\usepackage{amsfonts}
\usepackage{color}
\usepackage{theorem}

\usepackage{color}
\usepackage{ifthen}

\definecolor{red}{rgb}{.6,0,0}
\definecolor{green}{rgb}{0,.6,0}
\definecolor{darkgreen}{rgb}{0,0.3,0}

\newboolean{draft}

\newcommand{\cmt}[1]
{\ifthenelse {\boolean{draft}}
{{\sc \tiny \color{red} #1}}
{}}

\newcommand{\margincmt}[1]
{\ifthenelse {\boolean{draft}}
{\marginpar{{\sc \tiny \color{red} #1}}}
{}}
\newcommand{\inred}[1]
{\ifthenelse{\boolean{draft}}{{\color{red} #1}}{#1}}

\newcommand{\new}[1]
{\ifthenelse {\boolean{draft}}
{{\color{green} #1}}
{#1}}

\newcommand{\newb}[1]
{\ifthenelse {\boolean{draft}}
{{\color{blue} #1}}
{#1}}

\newcommand{\del}[1]
{\ifthenelse {\boolean{draft}}
{{\color{magenta} #1}}
{}}

\newboolean{details_on}

\newcommand{\details}[1]
{\ifthenelse {\boolean{details_on}}
{{\color{darkgreen} \tiny #1}}
{}}

\newcommand {\Map} {\mathbb{R}\mathbf{Map}}
\newcommand {\Parf} {\mathbb{R}\mathbf{Perf}}
\newcommand {\rh} {\mathbb{R}\underline{Hom}}
\newcommand {\rch} {\mathbb{R}\underline{\mathcal{H}om}}
\newcommand {\OO} {\mathcal{O}}
\newcommand {\Symp} {\mathbb{S}ymp}
\newcommand {\A} {\mathcal{A}}
\newcommand {\Spec} {\mathbf{Spec}}
\newcommand  {\dg}     {\mathbf{dg}}
\newcommand  {\mdg}     {\mathbf{dg}^{gr}}
\newcommand  {\edg}     {\epsilon-\mathbf{dg}}
\newcommand  {\medg}     {\epsilon-\mathbf{dg}^{gr}}

\newcommand  {\cdga}     {\mathbf{cdga}}

\newcommand  {\dAff}     {\mathbf{dAff}}
\newcommand  {\ncdga}     {\mathbf{cdga}^{\leq 0}}

\newcommand  {\ho}   {\mathrm{Ho}}

\newcommand  {\dSt}   {\mathbf{dSt}}

\newcommand{\s}{\infty}
\DeclareMathOperator*{\Holim}{\operatorname{{\sf Holim}}}

\newtheorem{thm}{Theorem}[section]
\newtheorem{prop}[thm]{Proposition}
\newtheorem{lem}[thm]{Lemma}
\newtheorem{df}[thm]{Definition}
\newtheorem{cor}[thm]{Corollary}

{\theorembodyfont{\rmfamily} \newtheorem{rmk}[thm]{Remark}}

\setboolean{draft}{true}
\setboolean{details_on}{true}

\begin{document}

\title{\textbf{
Shifted Symplectic Structures}}  
\author{\bf{T. Pantev}\footnote{Partially supported by NSF RTG grant
    DMS-0636606 and NSF 
grants DMS-0700446 and DMS-1001693.}
\footnote{Partially
    supported by the ANR grant 
ANR-09-BLAN-0151 (HODAG).} \\ \small{Department of Mathematics} \\
  \small{University of Pennsylvania}\\  
\and \bf{B. To\"en}\footnotemark[2] \\ \small{I3M} \\  \small{Universit\'e de
  Montpellier2}\\ \and \bf{M. Vaqui\'e}\footnotemark[2] \\
    \small{Institut de Math. de 
  Toulouse} \\  \small{Universit\'e Paul Sabati\'er}\\ \and
  \bf{G. Vezzosi}\footnotemark[2] \footnote{G. V. would like to thank
    the Insitute des Hautes \'Etudes Scientifiques (France) for
    providing perfect conditions for completing this paper.}
  \\ \small{Institut Math\'ematique de Jussieu} \\ 
  \small{Universit\'e Paris Diderot} }  

\date{March 2012}

\maketitle

\begin{abstract}
This is the first of a series of papers about \emph{quantization} in
the context of \emph{derived algebraic geometry}. In this first part,
we introduce the notion of \emph{$n$-shifted symplectic structures}
($n$-symplectic structures for short), a generalization of the notion
of symplectic structures on smooth varieties and schemes, meaningful
in the setting of derived Artin $n$-stacks (see \cite{hagII,seat}).
We prove that classifying stacks of reductive groups, as well as the
derived stack of perfect complexes, carry canonical $2$-symplectic
structures.  Our main existence theorem states that for any derived
Artin stack $F$ equipped with an $n$-symplectic structure, the derived
mapping stack $\textbf{Map}(X,F)$ is equipped with a canonical
$(n-d)$-symplectic structure as soon a $X$ satisfies a Calabi-Yau
condition in dimension $d$.  These two results imply the existence of
many examples of derived moduli stacks equipped with $n$-symplectic
structures, such as the derived moduli of perfect complexes on
Calabi-Yau varieties, or the derived moduli stack of perfect complexes
of local systems on a compact and oriented topological manifold. We
explain how the known symplectic structures on smooth moduli spaces of
simple objects (e.g. simple sheaves on Calabi-Yau surfaces, or simple
representations of $\pi_{1}$ of compact Riemann surfaces) can be
recovered from our results, and that they extend canonically as
$0$-symplectic structures outside of the smooth locus of simple
objects.  We also deduce new existence statements, such as the
existence of a natural $(-1)$-symplectic structure (whose formal
counterpart has been previously constructed in \cite{co,cg}) on the
derived mapping scheme $\textbf{Map}(E,T^*X)$, for $E$ an elliptic
curve and $T^*X$ is the total space of the cotangent bundle of a
smooth scheme $X$. Canonical $(-1)$-symplectic structures are also
shown to exist on Lagrangian intersections, on moduli of sheaves (or
complexes of sheaves) on Calabi-Yau $3$-folds, and on moduli of
representations of $\pi_{1}$ of compact topological
$3$-manifolds. More generally, the moduli sheaves on higher
dimensional varieties are shown to carry canonical shifted symplectic
structures (with a shift depending on the dimension).
\end{abstract}

\tableofcontents

\section*{Introduction} 
\addcontentsline{toc}{section}{Introduction}

This is the first part of a series of papers about \emph{quantization}
in the context of \emph{derived algebraic geometry}, and specifically
about the construction of quantized versions of various kinds of
moduli spaces.  In this article we start with the study of symplectic
structures in the derived setting by introducing the notion of
\emph{shifted symplectic structures of degree $n$} (or \emph{n-shifted
  symplectic structures}) where $n\in \mathbb{Z}$ is an arbitrary
integer. This is a direct and far reaching generalization of the
notion of symplectic structures on smooth algebraic varieties and
schemes (recovered when $n=0$), to the setting of derived and higher
derived Artin stacks of \cite{hagII,seat}.  In this work we give a
careful rigorous definition of $n$-shifted symplectic structures on
derived Artin stacks (see Definition \ref{d6}), and prove three
\emph{existence theorems} (see theorems \ref{t1}, \ref{t2}, \ref{t3})
which provide powerful construction methods and many examples. This
notion is an extension of the usual notion of symplectic structures on
smooth schemes on the one hand to higher algebraic stacks and on the
other hand to derived schemes and derived stacks.  Based on these
results, we recover some known constructions, such as the symplectic
structures on various types of moduli spaces of sheaves on surfaces
(see for instance \cite{mu,hule,in,iis}) and the symmetric obstruction
theories on moduli of sheaves on Calabi-Yau 3-folds (see
\cite{befan}), and prove new existence results, by constructing
natural $n$-shifted symplectic structures on many other moduli spaces,
including sheaves on higher dimensional varieties. These results may
be summarized as follows.

\begin{thm}\label{ti0}
\begin{enumerate}
\item Let $X$ be a smooth and proper Calabi-Yau variety of dimension
  $d$. Then the derived moduli stack of perfect complexes of
  quasi-coherent sheaves on $X$ admits a canonical $(2-d)$-shifted
  symplectic structure.
\item Let $X$ be a smooth and proper variety of dimension $d$. Then,
  the derived moduli stack of perfect complexes with flat
  connections on $X$ admits a canonical $2(1-d)$-shifted symplectic
  structure.
\item Let $M$ be a compact oriented topological manifold of dimension
$d$. Then, the derived moduli stack of perfect complexes of local
systems on $M$ admits a canonical $(2-d)$-shifted symplectic structure.
\end{enumerate}
\end{thm}

Future parts of this work will be concerned with the dual notion of
Poisson (and \emph{$n$-Poisson}) structures in derived algebraic
geometry, formality (and \emph{$n$-formality}) theorems, and finally
with quantization. 

\bigskip
\bigskip

\begin{center} 
\textbf{$p$-Forms, closed $p$-forms and symplectic forms in the derived
  setting} 
\end{center} 

A symplectic form on a smooth scheme $X$ (over some base ring $k$, of
characteristic zero), is the datum of a closed $2$-form $\omega \in
H^{0}(X,\Omega_{X/k}^{2,cl})$, which is moreover required to be
non-degenerate, i.e. it induces an isomorphism $\Theta_{\omega} :
T_{X/k} \simeq \Omega_{X/k}^{1}$ between the tangent and
cotangent bundles.  In our context $X$ will no longer be a scheme, but
rather a derived Artin stack in the sense of \cite{hagII,seat}, the
typical example being an $X$ that is the solution to some derived
moduli problem (e.g. of sheaves, or complexes of sheaves on smooth and
proper schemes, see \cite[Corollary 3.31]{tv}, or of maps between proper
schemes as in \cite[Corollary 2.2.6.14]{hagII}).  In this context,
differential $1$-forms are naturally sections in a quasi-coherent
complex $\mathbb{L}_{X/k}$, called the cotangent complex (see
\cite{il,seat}), and the quasi-coherent complex of $p$-forms is
defined to be $\wedge^{p}\mathbb{L}_{X/k}$. The $p$-forms on $X$ are
then naturally defined as sections of $\wedge^{p}\mathbb{L}_{X/k}$,
i.e. the set of $p$-forms on $X$ is defined to be the
(hyper)cohomology group $H^{0}(X,\wedge^{p}\mathbb{L}_{X/k})$. More
generally, elements in $H^{n}(X,\wedge^{p}\mathbb{L}_{X/k})$ are
called \emph{$p$-forms of degree $n$ on $X$} (see Definition \ref{d5}
and Proposition \ref{p3}).  The first main difficulty is to define the
notion of \emph{closed $p$-forms} and of \emph{closed $p$-forms of
degree $n$} in a meaningful manner. The key idea of this work is to
interpret $p$-forms, i.e. sections of $\wedge^{p}\mathbb{L}_{X/k}$, as
functions on the derived loop stack $\mathcal{L}X$ of \cite{seat,
chern, bena} by means of the HKR theorem of \cite{tove} (see also
\cite{bena}), and to interpret closedness as the condition of being
$S^{1}$-equivariant. One important aspect here is that $S^{1}$-equivariance
must be understood in the sense of homotopy theory, and therefore the 
closedness defined in this manner is not simply a property of a $p$-form but 
consists of an extra structure (see Definition \ref{d4}).
 This picture is accurate
(see Remark \ref{rloop1} and \ref{rloop2}), but technically difficult
to work with\footnote{One of the difficulties lies in the fact that we
need to consider only functions formally supported around the constant
loops $X \hookrightarrow \mathcal{L}X$, and this causes troubles
because of the various completions involved.}. We have therefore
chosen a different presentation, by introducing local constructions
for affine derived schemes, that are then glued over $X$ to obtain
global definitions for any derived Artin stack $X$.  With each
commutative dg-algebra $A$ over $k$, we associate a graded complex,
called the \emph{weighted negative cyclic complex of $A$ over $k$},
explicitly constructed using the derived de Rham complex of
$A$. Elements of weight $p$ and of degree $n-p$ of this complex are by
definition closed $p$-forms of degree $n$ on $\Spec\, A$ (Definition
\ref{d3}). For a general derived Artin stack $X$, closed $p$-forms are
defined by smooth descent techniques (Definition \ref{d5}). This
definition of closed $p$-forms has a more explicit local nature, but
can be shown to coincide with the original idea of $S^{1}$-equivariant
functions on the loop stack $\mathcal{L}X$ (using, for instance,
results from \cite{tove,bena}).

By definition a closed $p$-form $\omega$ of degree $n$ on $X$ has an
underlying $p$-form of degree $n$ (as we already mentioned this
underlying $p$-form does not determine the closed $p$-form $\omega$,
and several different closed $p$-form can have the same underlying
$p$-form). When $p=2$ this underlying $2$-form is an element in
$H^{n}(X,\wedge^{2}\mathbb{L}_{X/k})$, and defines a natural morphism
in the derived category of quasi-coherent complexes on $X$
$$\Theta_{\omega} : \mathbb{T}_{X/k} \longrightarrow \mathbb{L}_{X/k}[n],$$
where $\mathbb{T}_{X/k}$ is the tangent complex (i.e. the dual of
$\mathbb{L}_{X/k}$). With this notation in place, we give the main
definition of this paper:

\begin{df}\label{di1}
An \emph{$n$-shifted symplectic form} on a derived Artin stack $X$ is
a closed $2$-form $\omega$ of degree $n$ on $X$ such that the
corresponding morphism
$$
\Theta_{\omega} : \mathbb{T}_{X/k} \longrightarrow \mathbb{L}_{X/k}[n]
$$
is an isomorphism in the quasi-coherent derived category $D_{qcoh}(X)$
of $X$.  
\end{df}

When $n=0$ and $X$ is a \emph{smooth} $k$-scheme, the definition above
recovers the usual notion of a symplectic
structure, and nothing more.  Smooth
schemes do not admit $n$-shifted symplectic structures for $n\neq 0$, but
there are many interesting examples of $0$-shifted symplectic structures on
derived Artin stacks (see corollaries \ref{ct1}, \ref{ct3}; see also \cite{pecha}).
Therefore, not only the above definition provides an extension of the
notion of symplectic structure by introducing the parameter $n$, but
even for $n=0$ the notion of $0$-shifted symplectic structure is a new way to
extend the notion of symplectic structures on non-smooth
schemes. \\

An $n$-shifted symplectic form $\omega$ can be thought of as the data consisting
of a quasi-isomorphism 
$$
\Theta : \mathbb{T}_{X/k} \longrightarrow \mathbb{L}_{X/k}[n],
$$
together with an entire hierarchy of higher coherences expressing some
subtle relations between $\Theta$ and the differential geometry of
$X$. The quasi-isomorphism $\Theta$ can itself be understood as a kind
of duality between \emph{the stacky part of $X$}, expressed in the
non-negative part of $\mathbb{T}_{X}$, and \emph{the derived part of
$X$}, expressed in the non-positive part of $\mathbb{L}_{X/k}$ (this
is striking already when $n=0$, and this picture has to be qualified when $n$
is far away from $0$). In practice, when $X$ is some moduli of
sheaves on some space $M$, this duality is often induced by a version
of Poincar\'e duality (or Serre duality) on $M$, since tangent
complexes are then expressed in terms of the cohomology of $M$. It is
tempting to view $n$-shifted symplectic structures as a non-abelian
incarnation of Poincar\'e duality, which is definitely a good way to
think about them in the context of non-abelian cohomology (see the
paragraph \emph{Related works} at the end of this Introduction). 

\bigskip

\begin{center} \textbf{Existence results} \end{center}

In this paper we prove \emph{three} existence results for
$n$-shifted symplectic structures. To start with, we show that classifying
stacks $BG$, for reductive affine group schemes $G$, are naturally
endowed with $2$-shifted symplectic structures. The underlying $2$-form here
is clear, it is given by the degree $2$ shift of a non-degenerate
$G$-invariant 
quadratic form:
$$
\mathfrak{g}[1] \wedge \mathfrak{g}[1] \simeq Sym^{2}(\mathfrak{g})[2]
\longrightarrow k[2],
$$ where $\mathfrak{g}$ is the Lie algebra of $G$ over $k$. The fact
that this $2$-form can be naturally promoted to a closed $2$-form
on $BG$ follows from the simple observation that all $2$-forms on $BG$ are
canonically closed - geometrically this is due to the fact that the
space of functions on $\mathcal{L}BG$ is discrete, and thus $S^{1}$
acts on it in a canonically trivial manner. 

Our first existence result is the following

\begin{thm}\label{ti1}
The derived stack $\Parf$  of perfect complexes of quasi-coherent
sheaves is equipped with a natural 
$2$-shifted symplectic form. 
\end{thm}

The relation between the above theorem and the $2$-shifted symplectic form on
$BG$ is given by the canonical open embedding $B\mathrm{GL}_{n}
\subset \Parf$, sending a vector bundle to the corresponding perfect
complex concentrated in degree $0$: the $2$-shifted symplectic form on $\Parf$
restricts to the one on $B\mathrm{GL}_n$. The proof of Theorem
\ref{ti1} uses the Chern character for perfect complexes, with values
in negative cyclic homology. The weight $2$ part of the Chern
character of the universal perfect complex on $\Parf$ provides a
canonical $2$-form of degree $2$, which is non-degenerate by
inspection.

The second existence theorem we prove in this paper is a transfer of
$n$-shifted symplectic structures on a given derived Artin stack $F$ to the
derived mapping stack $\textbf{Map}(X,F)$, under certain
\emph{orientability} condition on $X$. This
statement can be viewed as an algebraic version of the 
AKSZ-formalism (see \cite{aksz} for the original reference), further
extended to the setting of derived Artin stacks.

\begin{thm}\label{ti2}
Let $X$ be a derived stack endowed with an $\OO$-orientation of
dimension $d$, and let $(F,\omega)$ be a derived Artin stack with an
$n$-shifted symplectic structure $\omega$. Then the derived mapping stack
$\textbf{Map}(X,F)$ carries a natural $(n-d)$-shifted symplectic structure.
\end{thm} 

The condition of having an $\OO$-orientation of dimension $d$ (see
Definition \ref{d9}) essentially means that $D_{qcoh}(X)$ satisfies
the Calabi-Yau condition in dimension $d$.  The typical
example is of course when $X$ is a smooth and proper Calabi-Yau scheme
(or Deligne-Mumford stack) of dimension $d$ (relative to $\Spec\,
k$). Other interesting examples are given, for instance, by de Rham or
Dolbeault homotopy types ($Y_{DR}$, $Y_{Dol}$ in the notation of
\cite{si1}) of a smooth and proper scheme $Y$ over $k$, for which
$\textbf{Map}(Y_{DR},F)$, or $\textbf{Map}(Y_{Dol},F)$, should be
understood as maps with flat connections, or with Higgs fields, from
$Y$ to $F$. Theorem \ref{ti1} and \ref{ti2} provide many examples of
$n$-shifted symplectic forms on moduli spaces of perfect complexes on
Calabi-Yau schemes, or flat perfect complexes, or perfect complexes
with Higgs fields, etc. The proof of Theorem \ref{ti2} is rather
natural, though the details require some care. We use the
evaluation morphism
$$
X \times \textbf{Map}(X,F) \longrightarrow F,
$$ 
to pull-back the $n$-shifted symplectic form on $F$ to a closed $2$-form of
degree $n$ on \linebreak $X \times \textbf{Map}(X,F)$. This closed $2$-form is
then \emph{integrated along} $X$, using the $\OO$-orientation (this is
a quasi-coherent integration, for which we need Serre duality), in
order to get a closed $2$-form of degree $(n-d)$ on
$\textbf{Map}(X,F)$. Then, we observe that this last $2$-form is
non-degenerate. A large part of the argument consists of defining
properly the integration map (see Definition \ref{d8}).

Finally, our third and last existence statement concerns symplectic
intersections and symplectic forms induced on them.  For this we
introduce the notion of a \emph{Lagrangian structure} on a morphism $L
\longrightarrow X$, where $X$ is equipped with an $n$-shifted symplectic form;
this is a generalization of the notion of Lagrangian submanifolds (a
closed immersion $L \hookrightarrow X$ of smooth schemes possesses a
Lagrangian structure if and only if $L$ is Lagrangian in $X$ in the
usual sense, and moreover this structure, if it exists, is unique).

\begin{thm}\label{ti3}
Let $(X,\omega)$ be a derived Artin stack with an $n$-shifted symplectic form
$\omega$, and let
$$
L \longrightarrow X \, ,\qquad L' \longrightarrow X
$$ 
be two morphisms of derived Artin stacks endowed with Lagrangian
structures. Then, the derived fiber product $L\times_{X}^{h}L'$
carries a natural $(n-1)$-shifted symplectic form.
\end{thm}

As a corollary, we see that the derived intersection of two Lagrangian
smooth subschemes $L,L' \subset X$, into a symplectic smooth scheme
$X$, always carries a natural $(-1)$-shifted symplectic structure.  

Before going further we would like to mention here that the use of
derived stacks in Theorems \ref{ti1}, \ref{ti2}, \ref{ti3} is crucial,
and that the corresponding results do not hold in the underived
setting. The reason for this is that if a derived Artin stack $F$ is
endowed with an $n$-shifted symplectic form $\omega$, then the pull-back of
$\omega$ to the truncation $h^{0}(F)$ is a closed $2$-form of degree
$n$ which is, in general, highly degenerate.

\bigskip

\begin{center} \textbf{Examples and applications} \end{center} 

The three theorems \ref{ti1}, \ref{ti2} and \ref{ti3} listed above,
imply the existence of many interesting and geometrically relevant
examples of $n$-shifted symplectic structures. For instance, let $Y$ be a
smooth and proper Deligne-Mumford stack with connected geometric
fibers of relative dimension $d$.

\begin{enumerate}
\item The choice of a fundamental class $[Y] \in
  H^{2d}_{DR}(Y,\OO)$\footnote{This stands for de Rham cohomology
  of $Y$ with coefficients in the trivial flat bundle $\OO_{Y}$. It can be computed  as usual, by 
  taking
  hypercohomology of $Y$ with coefficients in the de Rham complex.} in de Rham
  cohomology (relative to $\Spec\, k$) determines a canonical
  $2(1-d)$-shifted symplectic form on the derived stack
$$
\Parf_{DR}(Y):=\mathbf{Map}(Y_{DR},\Parf)
$$ 
of perfect complexes with flat connections on $Y$. 
\item The choice of a fundamental class $[Y] \in H^{2d}_{Dol}(Y,\OO)$\footnote{This stands for 
Dolbeault cohomology
  of $Y$ with coefficients in the trivial Higgs bundle $\OO_{Y}$. It can be computed  as usual, by 
  taking
  hypercohomology of $Y$ with coefficients in the Dolbeault complex.}
in Dolbeault cohomology (relative to $\Spec\, k$) determines a canonical
$2(1-d)$-shifted symplectic form on the derived stack
$$
\Parf_{Dol}(Y):=\mathbf{Map}(Y_{Dol},\Parf)
$$
of perfect complexes with Higgs fields.
\item The choice of a trivialization (when it exists)
$\omega_{Y/k} \simeq \OO_{Y}$,  
determines a canonical $(2-d)$-shifted symplectic form on the derived stack 
$$
\Parf(Y):=\mathbf{Map}(Y,\Parf)
$$
of perfect complexes on $Y$. 
\item If $M$ is a compact, orientable topological manifold of
  dimension $d$, then a choice 
of a fundamental class $[M]\in H_{d}(M,k)$ determines a canonical 
$(2-d)$-shifted symplectic form on the derived stack 
$$
\Parf(M):=\mathbf{Map}(M,\Parf)
$$ 
of perfect complexes on $M$\footnote{A perfect complex on $M$ is by
definition a complex of sheaves of $k$-modules locally
  quasi-isomorphic to a constant and bounded complex of sheaves made
  of projective $k$-modules of finite type.}.
\end{enumerate}

We note here that the derived stack of perfect complexes of
quasi-coherent sheaves on $Y$ considered above contains interesting
open substacks, such as the stack of vector bundles, or the stack of
simple objects. We use this observation, and our existence theorems,
to recover in a new and uniform way, some well known symplectic
structures on smooth moduli spaces of simple vector bundles (see
\cite{mu,in}), and on character varieties (see \cite{go,je}). One
corollary of our results states that these known symplectic structures
in fact \emph{extend} to $0$-shifted symplectic structures on the
ambient derived Artin stacks, and this explains what is happening to
the symplectic structures at the boundaries of these smooth open
substacks, i.e. at \emph{bad points} (vector bundles or
representations with many automorphisms, or with non-trivial
obstruction map, etc.).

As another application of our existence results we present a
construction of \emph{symmetric obstruction theories}, in the sense of
\cite{befan}, by showing that a $(-1)$-shifted symplectic structure on
a derived Artin stack $X$ always endows the truncation $h^{0}(X)$ with
a natural symmetric obstruction theory. This enables us to construct
symmetric obstruction theories on the moduli stack of local systems on
a compact topological $3$-manifold, or on the moduli stack of simple
perfect complexes on a Calabi-Yau $3$-fold. The latter result was
recently used by Brav-Bussi-Dupont-Joyce (\cite{joyce}) to prove that
the coarse moduli space of simple perfect complexes of coherent
sheaves, with fixed determinant, on a Calabi-Yau $3$-fold admits,
locally for the analytic topology, a \emph{potential}, i.e. it is
isomorphic to the critical locus of a function. This was a
longstanding problem in Donaldson-Thomas theory.  We remark that, as
shown recently by Pandharipande-Thomas (\cite{pandathomas}), such a
result is false for a general symmetric obstruction theory.  Hence the
existence of a local potential depends in a crucial way the existence
of a global $(-1)$-shifted symplectic structure on the derived moduli
stack of simple perfect complexes on a Calabi-Yau $3$-fold. The
results in \cite{joyce} suggest that one should have general
\emph{formal}, \emph{local-analytic} and perhaps \emph{\'etale local}
versions of the \emph{Darboux theorem} for $(-1)$-shifted - and maybe
even for general $n$-shifted - symplectic forms. It will be
interesting to compare such formal $(-1)$-shifted Darboux theorem with
the formal potential defined in \cite[Section~3.3]{koso}. Local
structure theorems of this type also hint at the existence of 
Donaldson-Thomas theory for Calabi-Yau manifolds of higher
dimensions. This is a completely unexplored territory with 
the first case being the case of $4$-folds, where the corresponding
derived moduli space carries a $(-2)$-symplectic structure by Theorem
\ref{ti1}.

Another interesting
$(-1)$-shifted symplectic form whose existence follows from our Theorem
\ref{ti2} is the one obtained on $\textbf{Map}(E,T^{*}X)$, where $E$
is an elliptic curve over $\mathrm{Spec}\, k$, $X$ is a
smooth $k$-scheme, and $T^{*}X$ is the total space of the cotangent
bundle of $X$ (relative to $k$) equipped with its canonical symplectic
structure. At the formal completion level this $(-1)$-shifted symplectic form
on $\textbf{Map}(E,T^{*}X)$ was constructed and studied in
\cite{co}. A nice feature of our construction is that it produces this
symplectic structure directly as a \emph{global} form on the derived scheme
$\textbf{Map}(E,T^{*}X)$. Specifically we have

\begin{cor}\label{ci1}
Let $E$ be an oriented elliptic curve over $k$ with a fixed algebraic volume form, and let $X$ be a smooth
$k$-scheme.  Then the derived mapping scheme $\textbf{Map}(E,T^{*}X)$
carries a canonical $(-1)$-shifted symplectic structure.
\end{cor}

\bigskip

\begin{center} \textbf{Future works and open questions} \end{center}

In a sequel to this paper, we will study the dual notion of
$n$\emph{-Poisson structures} on derived Artin stacks. This dual
notion is technically more delicate to handle than the $n$-shifted
symplectic structures discussed in this paper. This is essentially due
to the fact that it is a much less local notion. Nevertheless we can
follow the same reasoning and extract the notion of an $n$-Poisson
structure from the geometry of derived loop stacks and higher loop
stacks $\mathbb{L}^{(n)}X:=\mathbf{Map}(S^{n},X)$. As usual
$n$-shifted symplectic structures give rise to $n$-Poisson
structures, and correspond precisely to the \emph{non-degenerate}
$n$-Poisson structures. Without going into technical details one can
say that an $n$-Poisson structure on a derived Artin stack $X$
consists of the data of a \emph{bivector of degree $-n$}
$$
\mathcal{P} \in H^{-n}(X,\wedge^{2}(\mathbb{T}_{X/k}[n])[-2n])\simeq
H^{-n}(X,\phi_{n}^{2}(\mathbb{T}_{X/k})),
$$
where 
$$
\phi_{n}^{2}(\mathbb{T}_{X/k}) \simeq \left\{ 
\begin{array}{cc}
\wedge^{2}(\mathbb{T}_{X/k}) & if \; n\; is \; even \\
Sym^{2}(\mathbb{T}_{X/k}) & if \; n\; is \; odd
\end{array}
\right. ,
$$
together with higher coherences expressing a closedness
conditions. Now, it turns out that the 
complex
$$
\mathbb{R}\Gamma(X,Sym(\mathbb{T}_{X/k}[-n-1]))[n+2]
$$ can be identified with the tangent complex of the deformation
problem for the dg-category $L_{qcoh}(X)$ of quasi-coherent sheaves
considered as an \emph{$n$-fold monoidal dg-category} \footnote{We can
  make sense of this also for negative $n$, by shifting the formal
  deformation variable in degree $-2n$. For the moment, we will assume
  $n\geq 0$ in order to simplify the presentation.} : $n=0$
simply means as a dg-category, $n=1$ as a monoidal dg-category, $n=2$
as a braided monoidal dg-category, and so on. A higher version of
Kontsevich \emph{formality theorem} implies that an $n$-Poisson
structure $\mathcal{P}$ defines an element in this tangent complex,
satisfying a homotopy version of the master equation, and thus a
formal deformation of $L_{qcoh}(X)$ considered as an $n$-fold monoidal
dg-category. This formal deformation is, by definition, the
\emph{quantization of $X$ with respect to the $n$-Poisson structure
  $\mathcal{P}$}. As a consequence, if $X$ is endowed with an
$n$-shifted symplectic structure, then it has a canonical
quantization, which is, by definition, a formal deformation of
$L_{qcoh}(X)$ as an $n$-fold monoidal dg-category. This will be our
approach to construct quantizations of the derived moduli stacks on
which, in this paper, we have constructed $n$-shifted symplectic
structures. At the moment, this program is very much in progress, and
two main difficulties remain.  First of all the deformation theory of
dg-categories and $n$-fold monoidal dg-categories has not been fully
worked out in the literature. Even the case $n=0$ is not fully
understood, as explained in \cite{kelo} and \cite[Remark
  5.3.38]{lu}. We hope to make progress on this deformation theory by
introducing the new machinery of \emph{tame homological algebra}, and
\emph{tame dg-categories}. Ongoing work (see the forthcoming
\cite{tvv}) in that direction seems to provide a complete answer for
$n=0$. We hope that similar ideas will also work for arbitrary
$n$. The second issue concerns the higher version of Kontsevich
formality theorem mentioned above, which states that the natural Lie
bracket on the complex
$\mathbb{R}\Gamma(X,Sym(\mathbb{T}_{X/k}[-n-1]))[n+1]$ coming from the
fact that it is the tangent complex of the deformation problem of
$n$-fold monoidal dg-categories, is the natural one (i.e. equals an
appropriate version of the Schouten bracket). This higher version of
the formality theorem is, at present, still a conjecture, even for the
case $n=0$ and $X$ a general derived Artin stack (the only established
case is $n=0$ and $X$ a smooth scheme).

\bigskip

\begin{center} \textbf{Related works} \end{center} 

There are many related works that should be mentioned
here but for space reasons we will only discuss a small selection of
such works. 

To start with, our notion of $0$-shifted symplectic structure
generalizes the usual notion of symplectic structures on smooth
schemes to the setting of derived Artin stacks. In this context the
$0$-shifted symplectic structures provide a new point of view on
symplectic structures over non-smooth objects, and a comparison with
the notion of symplectic singularities and symplectic resolutions
\cite{kaledin.crepant,kls,namikawa11,nevins.mcgerty11} would certainly
be interesting. We note however that $0$-shifted symplectic structures
\emph{cannot} exist on singular schemes. This is due to the fact that
for the cotangent complex to be its own dual we need some non-trivial
stacky structure. Therefore, $0$-shifted symplectic structures do not
bring anything new for singularities of schemes, but are surely
interesting for singularities appearing on some moduli stacks and on
their coarse moduli spaces. As noted in \cite[Thm. 6.1]{fu}, there
exist coarse moduli spaces of sheaves on K3 surfaces having symplectic
singularities but with no symplectic resolution. However, these moduli
spaces are coarse moduli of natural derived Artin stacks, and
according to our results these derived stacks carry natural
$0$-shifted symplectic structures. It is tempting to consider the
derived Artin stack itself as a symplectic resolution of
its coarse moduli space.

We have explained how several well known symplectic structures can be
recovered from our existence theorem. Similar symplectic structures
are known to exist on moduli of certain sheaves on non Calabi-Yau
manifolds (see for instance \cite{kuma}). We believe that these can
also be recovered from a slight modification of our
constructions, but we will not pursue this direction in the
present work. 

Symplectic structures also 
appear in non-commutative geometry, in particular on moduli spaces
of sheaves on \emph{non-commutative Calabi-Yau varieties}. More generally, 
our existence Theorem \ref{t3} has an 
extension to the case of the derived stack $\mathcal{M}_{T}$ of objects in a Calabi-Yau dg-category $T$, 
constructed in\cite{tv}: the derived stack 
$\mathcal{M}_{T}$ carries a natural $(2-d)$-shifted symplectic structure when 
$T$ is a Calabi-Yau dg-category of dimension $d$. The proof of this non-commutative
extension of Theorem \ref{t3} is very close to the proof of our Theorem \ref{t1}, but 
is not included in this work. In another direction, 
we think that sheaves on non-commutative Calabi-Yau
varieties of fractional dimension should also carry
a suitable version of shifted symplectic structures.

Many of the $n$-shifted symplectic structures we construct in this paper live
on derived moduli stacks of bundles (or complexes of bundles) with
flat, or Higgs structures ($\mathbb{R}Loc_{DR}(X)$,
$\mathbb{R}Loc_{Dol}(X)$, etc.). These are the moduli stacks appearing
in \emph{non-abelian Hodge theory}, and in fact, the derived moduli of
flat perfect complexes $\Parf_{DR}(X)$ is used in \cite{si2} to
construct \emph{the} universal non-abelian Hodge filtration. The
$n$-shifted symplectic forms on these moduli stacks reflects the Poincar\'e
duality in de Rham (or Dolbeault, or Betti) cohomology, and this
should be thought of as an incarnation of Poincar\'e duality in
non-abelian cohomology.  They give important additional structures on
these moduli stacks, and are expected to play an important role for
the definition of \emph{polarizations} in non-abelian
Hodge theory.

Our Theorem \ref{ti2} should be viewed as an algebraic version of one
of the main construction in \cite{aksz}. Similarly, some of the
constructions and notions presented here are very close to
constructions and notions introduced in \cite{co}. More precisely, our
construction and results about degree $(-1)$ derived symplectic
structures might be seen as a \emph{globalization} of Costello's
\emph{formal derived} (in the sense of \cite{lu}) approach. A complete
comparison would require $\mathcal{C}^{\infty}$ and complex analytic
versions of derived algebraic geometry, and a notion of $n$-shifted
symplectic structures in such contexts. We are convinced that most, if
not all, the definitions and results we present in this work have
$\mathcal{C}^{\infty}$ and complex analytic analogs. Derived
differential and complex analytic geometries do exist thanks to
\cite{lu2}, but going through the notions of forms, closed forms and
symplectic structures in these settings is not completely
straightforward.

Finally, the reader will notice that most of the methods we give in
this work provide $m$-shifted symplectic structures starting from already
existing $n$-shifted symplectic structures with $n>m$ (e.g. Theorem \ref{t1}
and \ref{t2}). It would also be interesting to have constructions that
\emph{increase} the degree of symplectic structures. The
$2$-shifted symplectic forms on $BG$ and $\Parf$ are of this kind, but it
would be interesting to have general methods for constructing
$n$-shifted symplectic forms on quotients, dual to the one we have on fiber
products and on mapping stacks.

\bigskip

\bigskip

\noindent \textbf{Acknowledgements.} We are thankful to K.Costello
for sharing his ideas on $(-1)$-shifted symplectic forms with us, as
well as to D.Calaque and G.Ginot for conversations on
$E_{n}$-deformation theory and the higher dimensional formality
conjecture.  We are grateful to D.Joyce for sending us a draft of
\cite{joyce}, where $(-1)$-shifted symplectic forms were used to
establish the existence of a local analytic potential in
Donaldson-Thomas theory for Calabi-Yau $3$-folds; this was a great
happy ending to our paper. We also thank F.Bonechi, A.Cattaneo,
D.Kaledin, D.Kazhdan, M.Kontsevich, D.Maulik, F.Paugam and
B.Tsygan for stimulating conversations on the subject while we were
writing this paper.  We finally thank P-E.Paradan for pointing out
the reference \cite{je}.

\bigskip
\bigskip

\noindent \textbf{Notations and conventions.} 
\begin{itemize}

\item $k$ is a base commutative ring, noetherian and of residual
characteristic zero.

\item $dg_{k}$ is the category of dg-modules over $k$ (i.e. of
complexes of $k$-modules). By convention, the differential of an
object in $dg_{k}$ \emph{increases} degrees. For an object $E \in
dg_{k}$, we will sometimes use the notation
$$\pi_{i}(E):=H^{-i}(E).$$

\item $cdga_{k}$ is the category of commutative dg-algebras over $k$,
and $cdga_{k}^{\leq 0}$ its full subcategory of non-positively graded
commutative dg-algebras.

\item $dg_{k}$, $cdga_{k}$ (respectively $cdga_{k}^{\leq 0}$) are
  endowed with their natural model structures for which equivalences
  are quasi-isomorphisms, and fibrations are epimorphisms
  (respectively epimorphisms in strictly negative degrees).

\item $dAff_{k}:=(cdga_{k}^{\leq 0})^{op}$ is the category of derived
affine $k$-schemes.

\item The expression \emph{$\s$-category} will always refer to
\emph{$(\infty, 1)$-category} (see \cite{ber}). To fix ideas we will
use Segal categories as models for $\s$-categories (see \cite{si3},
and \cite[\S 1]{chern}).

\item The $\s$-categories associated to the model categories $dg_{k}$,
$cdga_{k}^{\leq 0}$, $dAff_{k}$ are denoted by $\dg_{k}$,
$\ncdga_{k}$, $\dAff_{k}$.

\item The $\s$-category of simplicial sets is denoted by
$\mathbb{S}$. It is also called the $\s$-category of spaces, and
\emph{space} will be used to mean \emph{simplicial set}.

\item The $\s$-category of derived stacks over $k$, for the \'etale
topology, is denoted by $\dSt_{k}$ (see \cite{seat,hagII}). If $X$ is
a derived stack, the $\s$-category of derived stacks over $X$ is
denoted by $\dSt_{X}$. The truncation of a derived stack $X$ is
denoted by $h^{0}(X)$. The derived mapping stack between $X$ and $Y$
is denoted by $\Map(X,Y)$.

\item We will use the expressions \emph{homotopy limits} and
  \emph{$\s$-limits} interchangeably to refer either to homotopy
  limits in an ambient model category, or to limits in an ambient
  $\s$-category. The same convention will be used for
  colimits. Homotopy fiber products will be denoted as usual by
  $X\times_{Z}^{h}Y$.

\item The mapping space between two objects $a$ and $b$ in an
$\s$-category $A$ is denoted by $Map_{A}(a,b)$. Points in
$Map_{A}(a,b)$ will be called \emph{morphisms in $A$}, and paths in
$Map_{A}(a,b)$ \emph{homotopies}. A morphism in an $\s$-category will
be called an \emph{equivalence} if it is a homotopy equivalence
(i.e. becomes an isomorphism in the homotopy category).  The word
\emph{equivalence} will also refer to a weak equivalence in a model
category.

\item A \emph{derived Artin stack} is by definition a derived stack
which is $m$-geometric for some integer $m$, for the \'etale topology
and the class $\mathbf{P}$ of smooth maps (see \cite[\S 2.2.2]{hagII}
and \cite{seat}). All derived Artin stacks are assumed to be locally
of finite presentation over $\mathrm{Spec}\, k$.

\item For a derived stack $X$, its quasi-coherent derived category is
denoted by $D_{qcoh}(X)$, and its $\s$-categorical version by
$L_{qcoh}(X)$. The larger derived category of all $\OO_{X}$-modules is
denoted by $D(\OO_{X})$, and its $\s$-categorical version by
$L(\OO_{X})$. By definition, $L_{qcoh}(X)$ is a full sub-$\s$-category
of $L(\OO_{X})$. As usual, morphisms between $x$ and $y$ in
$D_{qcoh}(X)$ will be denoted by $[x,y]$. We have
$$
[x,y]\simeq \pi_{0}(Map_{L_{qcoh}(X)}(x,y)).
$$ We refer to \cite{seat,deraz}, and also \cite{lu2}, for detailed
definitions and properties of these $\s$-categories.

\item Complexes of morphisms between objects $x$ and $y$ in
$L(\OO_{X})$ or $L_{qcoh}(X)$ will be denoted by $\rh(x,y)$. The
$\s$-categories $L(\OO_{X})$ and $L_{qcoh}(X)$ have natural symmetric
monoidal structures (in the sense of \cite{chern,lu4}), and this
monoidal structures are closed. The internal mapping object in
$L(\OO_{X})$ between
$x$ and $y$ is denoted by $\rch(x,y)$. In particular, for an object
$x$, we denote by $x^{\vee}:=\rch(x,\OO_{X})$ its dual. Perfect
complexes on $X$ are by definition dualizable objects in
$L_{qcoh}(X)$.

\item Expressions such as $f^{*}$, $f_{*}$ and $\otimes_{k}$ should be
understood in the derived sense, and should be read as
$\mathbb{L}f^{*}$, $\mathbb{R}f_{*}$, $\otimes_{k}^{\mathbb{L}}$,
unless specified otherwise.

\item For a derived Artin stack $X$, we denote by $\mathbb{L}_{X/k}
  \in L_{qcoh}(X)$ its cotangent complex over $\mathrm{Spec}\, k$ (see
  \cite{seat,hagII} and \cite[\S 7.3]{lu4}). Since, according to our
  conventions, $X$ is assumed to be locally of finite presentation, it
  follows that $\mathbb{L}_{X/k}$ is a perfect complex on $X$, and
  thus is dualizable. Its dual is denoted by
  $\mathbb{T}_{X/k}:=\mathbb{L}_{X/k}^{\vee}$ and is called the
  tangent complex of $X$ over $\mathrm{Spec}\, k$.
\end{itemize}

\section{Definitions and properties}

In this first section we give the definitions of \emph{$p$-forms},
\emph{closed $p$-forms} and \emph{symplectic forms} over a derived
Artin stack. We will start by some elementary constructions in the
setting of mixed and graded mixed complexes.  These constructions will
then be applied to define $p$-forms and closed $p$-forms over an
affine derived scheme $\Spec\, A$, by using some explicit graded mixed
complexes constructed from the derived de Rham complex of
$A$. Finally, these constructions are shown to be local for the smooth
topology and to glue on smooth covers, giving global notions for any
derived Artin stack. This smooth descent property is not a completely
obvious statement, and its proof requires some care.

\subsection{Graded mixed complexes}\label{gmc}

Recall from \cite{kas} that a \emph{mixed complex over $k$} is a
dg-module over $k[\epsilon]=H_{*}(S^1,k)$, where $\epsilon$
is in degree $-1$ and satisfies $\epsilon^{2}=0$. The category of
mixed complexes will be denoted by $\epsilon-dg_{k}$, and also called
the category of $\epsilon$-dg-modules over $k$. The differential of
the complex of $k$-modules underlying an object $E\in \epsilon-dg_{k}$
will be denoted by $d$, and by convention it raises degrees $d : E^{n}
\longrightarrow E^{n+1}$.

The tensor product $\otimes_{k}$ makes $\epsilon-dg_{k}$ into a
symmetric monoidal category, which is moreover a symmetric monoidal
model category for which the weak equivalences are the
quasi-isomorphisms (see \cite[\S 2]{tove}).  By definition a
\emph{graded mixed complex over $k$} is a mixed complex $E$ over $k$,
equipped with a direct sum decomposition of the underlying complex of
$k$-dg-modules
$$
E:=\bigoplus_{p \in \mathbb{Z}}E(p),
$$
in such a way that multiplication by $\epsilon$ has degree $1$ 
$$
\epsilon : E(p) \longrightarrow E(p+1),
$$
while the differential $d$ of the complex of $k$-modules underlying $E$ 
respects this grading
$$
d : E(p) \longrightarrow E(p).
$$ The extra grading will be called the \emph{weight grading}, and we
will say that elements in $E(p)$ have \emph{weight} $p$.

Graded mixed complexes over $k$ form a category denoted by
$\epsilon-dg^{gr}_{k}$, that is again a symmetric monoidal model
category, for the tensor product over $k$ with the usual induced
grading
$$
(E\otimes_{k}F)(p):=\bigoplus_{i+j=p}E(i)\otimes_{k}E(j),
$$
and with weak equivalences being the quasi-isomorphisms. 

\begin{rmk}\label{r0} Alternatively, graded mixed complexes can be
  viewed as
  dg-comodules over the commutative dg-Hopf-algebra (see also \cite{bena})
$$
B_{\epsilon}:=H^{*}(\mathbb{G}_{m} \ltimes
B\mathbb{G}_{a},\mathcal{O}),
$$ 
the cohomology Hopf algebra of the
semi-direct product group stack $\mathbb{G}_{m} \ltimes
B\mathbb{G}_{a}$. It is the semi-direct (or cross-)product of the
multiplicative Hopf algebra $k[t,t^{-1}]$ with $k[e]=H^{*}(S^1,k)$, via
the natural action of $\mathbb{G}_{m}$ on $k[e]$ given by rescaling
$e$. More precisely, as a commutative dg-algebra, $B_{\epsilon}$ is
$k[t,t^{-1}]\otimes_{k} k[e]$, with zero differential and 
comultiplication determined by
$$
\Delta(t)= t\otimes t \qquad
\Delta(e)=t\otimes e.
$$
Note that the tautological equivalence of symmetric monoidal categories
$$
\epsilon-dg^{gr}_{k} \simeq B_{\epsilon}-dg-comod,
$$
commutes with the two forgetful functors 
to the category of complexes of $k$-dg-modules. 
\end{rmk}

Let $E \in \epsilon-dg_{k}$ be a mixed complex over $k$. We may form
the usual negative cyclic object $NC(E)$, which is a dg-module over
$k$ defined in degree $n$ by the formula
$$
NC^{n}(E):=\prod_{i\geq 0} E^{n-2i}.
$$
For an element $\{m_{i}\}_{i} \in \prod_{i\geq 0} E^{n-2i}$, 
the differential is defined by the formula
$$
D(\{m_{i}\})_{j}=\epsilon m_{j-1}+dm_{j}.
$$
If $E$ is moreover a \emph{graded} mixed complex, then the extra
grading $E=\bigoplus_{p}E(p)$ defines 
a sequence of sub-complexes $NC(E)(p) \subset NC(E)$, by 
$$
NC^{n}(E)(p):=\prod_{i\geq 0} E^{n-2i}(p+i).
$$ 
These sub-complexes are natural direct summands (i.e the inclusions
admit natural retractions), and thus we have natural morphisms
$$
\bigoplus_{p}NC(E)(p) \longrightarrow NC(E) \longrightarrow
\prod_{p}NC(E)(p).
$$
These are monomorphisms of complexes, but in general they are neither
isomorphisms nor quasi-isomorphisms.

The above constructions are functorial in $E$, and define a functor
$$
NC : \epsilon-dg^{gr}_{k} \longrightarrow  dg_{k},
$$
from the category of graded mixed complexes to the 
category of complexes. Moreover, the functors $E \mapsto NC(E)(p)$ defines 
direct summands of the functor $NC$, together with natural transformations
$\bigoplus_{p}NC(p) \longrightarrow NC \longrightarrow \prod_{p}NC(p)$.
We will be interested in the family of functors $NC(p) :
\epsilon-dg^{gr}_{k} \longrightarrow dg_{k}$, 
as well as in their direct sum 
$$
NC^{w}:=\bigoplus_{p}NC(p) : \epsilon-dg^{gr}_{k} \longrightarrow
dg^{gr}_{k},
$$ 
where $dg^{gr}_{k}$ denotes the category of graded complexes of
$k$-modules - i.e. complexes of $k$-modules equipped with an extra
$\mathbb{Z}$-grading, with morphisms the maps of complexes preserving
the extra grading.  

\begin{df}\label{d1}
For a graded mixed complex $E$, the \emph{weighted negative cyclic
complex} is defined by
$$
NC^{w}(E):=\bigoplus_{p\in \mathbb{Z}}NC(E)(p) \in dg^{gr}_{k}.
$$
Its cohomology is called the \emph{weighted negative cyclic homology of
  $E$}, and is denoted by 
$$
NC^{w}_{n}(E):=H^{-n}(NC^{w}(E)) \qquad
NC^{w}_{n}(E)(p):=H^{-n}(NC^{w}(E)(p)).
$$
The natural decomposition
$$
NC^{w}_{n}(E)\simeq \bigoplus_{p\in \mathbb{Z}}NC^{w}_{n}(E)(p)
$$
is called the \emph{Hodge decomposition}.
\end{df}

As observed above, we have a natural morphism
$$NC^{w}(E) \longrightarrow NC(E)$$ which is, in general, not a
quasi-isomorphism. The complex $NC(E)$ computes the usual
\emph{negative cyclic homology} of the mixed complex $E$, which can
differ from the weighted cyclic homology as defined above. The two
coincide under obvious boundedness conditions on $E$.

Both categories $\epsilon-dg^{gr}_{k}$ and $dg^{gr}_{k}$ have natural
\emph{projective} model structures, where weak equivalences are
quasi-isomorphisms and fibrations are epimorphisms, on the underlying
complexes of $k$-modules (see \cite{tove}).

\begin{prop}\label{p1}
The functor
$$NC^{w} : \epsilon-dg^{gr}_{k} \longrightarrow dg^{gr}_{k}$$
is a right Quillen functor for the projective model structures.
\end{prop}
\textbf{Proof.} We denote by $k(p) \in \epsilon-dg^{gr}_{k}$ the graded
mixed complex consisting of $k$ sitting in degree $0$ and having pure
weight $p$. We construct an explicit cofibrant replacement $Q(p)$ of
$k(p)$ as follows. As a graded
$k[\epsilon]$-module, $Q(p)$ is
$$
Q(p)=\bigoplus_{j \geq
  0}k[\epsilon][2j]=k[\epsilon][a_{0},\dots,a_{j},\dots],
$$ 
where $a_{j}$ stands for the canonical generator of $k[\epsilon][2j]$
(in degree $-2j$).  
The differential in $Q(p)$ is then given by 
$$
da_{j}+\epsilon\cdot a_{j-1}=0.
$$
Finally, $a_{j}$ is declared to be of weight $p+j$. The natural 
projection
$$
Q(p) \longrightarrow k(p)
$$
sending all the $a_j$'s to zero for $j>0$, and $a_0$ to $1 \in k$, is
a quasi-isomorphism of graded mixed complexes. Moreover $Q(p)$ is a
free graded module over $k[\epsilon]$, and thus is a cofibrant object
in $\epsilon-dg^{gr}_{k}$. Finally, we have
$$
NC^{w}(E)(p)\simeq \underline{Hom}(Q(p),E),
$$ 
where $\underline{Hom}$ stands for the complex of morphisms between
two objects in $\epsilon-dg^{gr}_{k}$ ($\epsilon-dg_{k}^{gr}$ is a
$C(k)$-model category in a natural way, where $C(k)$ is the symmetric
monoidal model category of complexes of $k$-modules). This finishes
the proof of the proposition.  \hfill $\Box$  \\

Since all objects in the model categories $\epsilon-dg^{gr}_{k}$ and
$dg^{gr}_{k}$ are fibrant for the projective model structures,
Proposition \ref{p1} implies in particular that $NC^{w}$ preserves
quasi-isomorphisms. Following our conventions, we will denote by
$$
\mdg_{k}, \qquad \edg_{k}, \qquad \text{and\,\,\,\,\,\,} \medg_{k},
$$ 
the $\s$-categories obtained by Dwyer-Kan localization from the
corresponding model categories $dg^{gr}_{k}$, $\epsilon-dg_{k}$ and
$\epsilon-dg^{gr}_{k}$.  The functor $NC^{w}$ thus defines an
$\s$-functor
$$
NC^{w}=\bigoplus_{p}NC^{w}(p) : \medg_{k} \longrightarrow \mdg_{k}.
$$

\begin{cor}\label{cp1}
The $\s$-functor
$$
NC^{w} : \medg_{k} \longrightarrow \mdg_{k}
$$
preserves $\s$-limits. Moreover, we have natural equivalences
$$
NC^{w}(E)(p)\simeq \mathbb{R}\underline{Hom}(k(p),E)
$$ 
where $k(p) \in \epsilon-dg^{gr}_{k}$ is the graded mixed complex
consisting of $k$ sitting in degree $0$ and having pure weight $p$.
\end{cor}
\textbf{Proof.} This is a consequence of the statement and proof of
proposition \ref{p1}. \hfill $\Box$  

\

\smallskip

\begin{rmk}\label{r1}
The notion of weighted negative cyclic homology has the following
geometric interpretation (see also \cite{bena}). We let the
multiplicative group scheme $\mathbb{G}_{m}$ (over $\Spec\, k$) act on
the group stack $B\mathbb{G}_{a}$, and form the semi-direct product
group stack $\mathbb{G}_{m}\ltimes B\mathbb{G}_{a}$. It can be shown
from the definitions that there exists an equivalence of
$\s$-categories
$$
L_{qcoh}(B(\mathbb{G}_{m}\ltimes B\mathbb{G}_{a})) \simeq \medg_{k}.
$$ 
This equivalence is obtained by sending a quasi-coherent complex
$E$ on the stack $B(\mathbb{G}_{m}\ltimes B\mathbb{G}_{a})$ to its
fiber at the base point $\Spec\, k \longrightarrow
B(\mathbb{G}_{m}\ltimes B\mathbb{G}_{a})$, which is naturally equipped
with an action of $\mathbb{G}_{m}\ltimes B\mathbb{G}_{a}$, and thus a
structure of a comodule over the dg-coalgebra of cochains
$C^{*}(\mathbb{G}_{m}\ltimes B\mathbb{G}_{a},\mathcal{O})$.  This
dg-coalgebra turns out to be formal and quasi-isomorphic to the
dg-coalgebra $B_{\epsilon}$ discussed above - the semi-direct product
of $k[t,t^{-1}]$ with $k[e]$. This fiber is therefore a graded mixed
complex, and it is straightforward to check that this construction
induces an 
equivalence of $\s$-categories as stated.  Using this point of view,
the $\s$-functor $NC^{w}$ has the following interpretation. Consider
the natural projection
$$ 
\pi : B(\mathbb{G}_{m}\ltimes B\mathbb{G}_{a}) \longrightarrow 
B\mathbb{G}_{m}.
$$
It induces a direct image on the $\s$-categories of quasi-coherent
complexes 
$$
\pi_{*} : L_{qcoh}(B(\mathbb{G}_{m}\ltimes B\mathbb{G}_{a}))
\longrightarrow  
L_{qcoh}(B\mathbb{G}_{m}), 
$$
right adjoint to the pull-back functor 
$$
\pi^{*} : L_{qcoh}(B\mathbb{G}_{m}) \longrightarrow
L_{qcoh}(B(\mathbb{G}_{m}\ltimes B\mathbb{G}_{a})).
$$ If we identify $L_{qcoh}(B(\mathbb{G}_{m}\ltimes B\mathbb{G}_{a}))$
with $\medg_{k}$ as above, and $L_{qcoh}(B\mathbb{G}_{m})$ with
$\mdg_{k}$, then $\pi_{*}$ becomes \emph{isomorphic} to $NC^{w}$. In
other words, $NC^{w}(E)$ computes the homotopy fixed point of $E$
under the action of $B\mathbb{G}_{a}$.  The residual
$\mathbb{G}_{m}$-action then corresponds to the grading by the pieces
$NC^{w}(E)(p)$.
\end{rmk}

\

\smallskip

We finish this section by observing that the graded complex
$NC^{w}(E)$ comes equipped with a natural projection
$$
NC^{w}(E) \longrightarrow E,
$$
which is functorial in $E$ and is a morphism of graded complexes. It
consists of  
the morphisms induced by the projection to the $i=0$ component
$$
NC^{w}(E)^{n}(p)=\prod_{i\geq 0} E^{n-2i}(p+i) \longrightarrow 
E^{n}(p).
$$

\begin{rmk}\label{epsilond} If $E$ is a graded $\epsilon$-dg-module
  (i.e a graded mixed complex) over $k$, then for any weight $p$ there
  is 
  a natural map - that might be called the absolute
  $\epsilon$-differential -  
$$
D_{E}(p): E(p) \longrightarrow
  \mathrm{NC}^{w}(E)(p+1)[-1]
$$ 
whose degree $m$ piece is 
$$
D_{E}^{m}(p):
  E^{m}(p) \rightarrow \prod_{i\geq 0} E^{m-2i-1}(p+1+i): x_{m,p}
  \mapsto (\epsilon x_{m,p}, 0,0 ,\ldots).
$$  
\end{rmk}

\subsection{$p$-Forms, closed $p$-forms and $n$-shifted symplectic
  structures} 

Let $A \in cdga^{\leq 0}_{k}$ be a commutative differential
non-positively graded algebra over $k$, and let us denote by
$\Omega^{1}_{A/k}$ the $A$-dg-module of K\"ahler differential
$1$-forms of $A$ over $k$.  We can form its de Rham algebra over $k$
(see \cite{tove})
$$
DR(A/k):=Sym_{A}(\Omega^{1}_{A/k}[1]),
$$ 
where in contrast with our usual usage $Sym_{A}$ here refers to the
underived symmetric product of $A$-dg-modules. The complex $DR(A/k)$
is in a natural way a commutative $\epsilon$-dg-algebra in the sense
of \cite[\S 2]{tove}. Here we will consider $DR(A/k)$ \emph{just} as a
graded mixed complex, by forgetting the extra multiplicative
structure. The underlying complex of $k$-modules is simply
$$
Sym_{A}(\Omega^{1}_{A/k}[1]) \simeq \bigoplus_{p} \Omega^{p}_{A/k}[p],
$$
where $\Omega^{p}_{A/k}:=\wedge^{p}_{A}\Omega^{1}_{A/k}$. The mixed
structure is induced by the de Rham differential $\epsilon:=d_{DR} :
\Omega^{p}_{A/k} \longrightarrow \Omega^{p+1}_{A/k}$, and is the
unique mixed structure on $Sym_{A}(\Omega^{1}_{A/k}[1])$ making it
into an $\epsilon$-cdga and for which the action of $\epsilon$ on the
factor $A=\wedge^{0}_{A}\Omega^{1}_{A/k}$ is the universal
dg-derivation $d : A \longrightarrow \Omega^{1}_{A/k}$ (see \cite[\S
2]{tove}).  Finally, the grading on $DR(A/k)$ is the one for which
$DR(A/k)(p):=\Omega^{p}_{A/k}[p]$.  This graded mixed structure on
$DR(A/k)$ is also compatible with the multiplicative structure, and
makes it into a graded mixed cdga over $k$, but we will not make use
of this finer structure in this work.

The assignment $A \mapsto DR(A/k)$ defines a functor
$$
cdga^{\leq 0}_{k} \longrightarrow \epsilon-dg^{gr}_{k}.
$$
This functor can be derived on the left, by pre-composing it with
a cofibrant replacement functor on $cgda^{\leq 0}_{k}$, to obtain
$$
\mathbb{L}DR(-/k) : cdga^{\leq 0}_{k} \longrightarrow
\epsilon-dg^{gr}_{k}
$$ 
which now preserves quasi-isomorphisms. Therefore, it induces a
well defined $\s$-functor $\mathbf{DR}(-/k) : \cdga_{k}
\longrightarrow \medg_{k}$.  Recall that, for $A$ a cdga over $k$, the
underlying complex of $\mathbf{DR}(A/k)$ is
$$\bigoplus_{p}(\wedge^{p}_{A}\mathbb{L}_{A/k})[p],$$
where $\mathbb{L}_{A/k}$ is the cotangent complex of $A$ over $k$, and 
$\wedge_{A}^{p}$ must now be understood in the derived sense (see
the proof of \cite[Proposition 2.4]{tove}).

\begin{df}\label{d2}
Let $A\in \cdga_{k}$. The \emph{weighted negative cyclic complex of
  $A$ over $k$} is defined by 
$$
NC^{w}(A/k):=NC^{w}(\mathbf{DR}(A/k)).
$$
This defines an $\s$-functor
$$
NC^{w} : \cdga_{k} \longrightarrow \mdg_{k}.
$$
\end{df}

As we have seen, for any graded mixed complex $E$, we have a natural
projection
$$
NC^{w}(E) \longrightarrow E,
$$
which is a morphism of graded complexes. We get this way, for any
$p\geq 0 $, a natural morphism of complexes 
$$
NC^{w}(A/k)(p) \longrightarrow  \wedge_{A}^{p}\mathbb{L}_{A/k}[p].
$$

For a complex of $k$-modules $E$, we will denote by $|E|$ the
simplicial set obtained by the Dold-Kan correspondence (applied to the
truncation $\tau_{\leq 0}(E)$). By definition we have a natural weak
equivalence $|E| \simeq Map_{\dg_{k}}(k,E)$, where $Map_{\dg_{k}}$
denotes the mapping space (i.e. simplicial set) in the $\infty$-category 
$\dg_{k}$.  For $A\in \cdga_{k}$, and two integers $p\geq 0$
and $n\in \mathbb{Z}$, we set
$$
\A^{p}_{k}(A,n):=|\wedge^{p}_{A}\mathbb{L}_{A/k}[n]| \in \mathbb{S}.
$$
This defines an $\s$-functor $\A^{p}_{k}(-,n) : \cdga_{k}
\longrightarrow \mathbb{S}$. 
In the same way, we set
$$
\A^{p,cl}_{k}(A,n):=|NC^{w}(A/k)[n-p](p)|.
$$
Using the natural projection mentioned above
$NC^{w}(A/k)[n-p](p) \longrightarrow \wedge^{p}\mathbb{L}_{A/k}[n]$, 
we deduce a natural morphism
$\A^{p,cl}_{k}(A,n) \longrightarrow \A^{p}_{k}(A,n)$.
We have therefore two $\s$-functors
$$
\A^{p,cl}_{k}(-,n) \, , \,\, \A^{p}_{k}(-,n) : \cdga_{k}
\longrightarrow \mathbb{S},
$$
together with a natural transformation
$\A^{p,cl}_{k}(-,n) \longrightarrow \A^{p}_{k}(-,n)$.

\begin{df}\label{d3}
For $A\in \cdga_{k}$, 
the simplicial set $\A^{p}_{k}(A,n)$ (respectively
$\A^{p,cl}_{k}(A,n)$) is called
the \emph{space of $p$-forms of degree $n$ on the derived stack
  $\Spec\, A$, relative to $k$} 
(respectively the \emph{space of closed $p$-forms of degree $n$ on the
  derived 
  stack $\Spec\, A$, relative to $k$}).  
\end{df}

As usual, when the ground ring $k$ is clear from the context, we will
simply write $\A^{p}(-,n)$ and $\A^{p,cl}(-,n)$ for $\A_{k}^{p}(-,n)$
and $\A_{k}^{p,cl}(-,n)$.

\begin{rmk}\label{rloop1}
We have seen that graded mixed complexes can be understood as
quasi-coherent complexes on the stack $B\mathcal{H}$, where
$\mathcal{H}$ is the semi-direct product group stack
$\mathbb{G}_{m}\ltimes B\mathbb{G}_{a}$ (see our remark \ref{r1}). If
we continue this point of view, $p$-forms and closed $p$-forms can be
interpreted as follows (see also \cite{bena}).

Let $X=\Spec\, A$ be a derived affine scheme over $k$, and consider its 
derived loop stack $\mathcal{L}X:=\Map(S^{1},X)$. The natural
morphism $S^{1}=B\mathbb{Z} \longrightarrow B\mathbb{G}_{a}$, 
induced by the inclusion $\mathbb{Z} \hookrightarrow \mathbb{G}_{a}$,
induces a morphism 
$$
\mathcal{L}^{u}X:=\Map(B\mathbb{G}_{a},X) \longrightarrow
\mathcal{L}X.
$$
This morphism turns out to be an equivalence of derived
schemes. Therefore, the  
group stack $\mathcal{H}$ of automorphisms of $B\mathbb{G}_{a}$, acts naturally 
on $\mathcal{L}X$. We can form the quotient stack and consider the
natural projection 
$$
p : [\mathcal{L}X/\mathcal{H}] \longrightarrow B\mathcal{H}.
$$ 
Using the results of \cite{tove} it is possible to show that there
exists a functorial equivalence in $L_{qcoh}(B\mathcal{H})$
$$
p_{*}(\mathcal{O}_{[\mathcal{L}X/\mathcal{H}]}) \simeq
\mathbf{DR}(A/k),
$$ 
where $\mathbf{DR}(A/k)$ is viewed as an object in
$L_{qcoh}(B\mathcal{H})$ using our remark \ref{r1}.

As $p$-forms and closed $p$-forms are defined directly from the graded
mixed complexes, this explains the precise relation between our notion
of closed $p$-forms, and functions on derived loop stacks. For
instance, we have that $NC^{w}(A/k)$ are simply the
$B\mathbb{G}_{a}$-invariants (or equivalently the $S^{1}$-invariants
through $S^{1} \longrightarrow B\mathbb{G}_{a}$) in the complex of
functions on $\mathcal{L}X$, in the sense that we have a natural
equivalence of quasi-coherent sheaves on $B\mathbb{G}_{m}$
$$
q_{*}(\mathcal{O}_{[\mathcal{L}X/\mathcal{H}]}) \simeq NC^{w}(A/k),
$$
where now $q$ is the projection $[\mathcal{L}X/\mathcal{H}]
\longrightarrow  
B\mathcal{H} \longrightarrow B\mathbb{G}_{m}$. 
\end{rmk}

Note that the space $\A^{p,cl}(A,n)$, of closed $p$-forms of some
degree $n$, is not a full sub-space (i.e. not a union of connected
components) of $\A^{p}(A,n)$.  For a point $w\in \A^{p}(A,n)$, the
homotopy fiber $K(w)$ of the map $\A^{p,cl}(A,n) \longrightarrow
\A^{p}(A,n)$, taken at $w$, can be a complicated space. Contrary to
what the terminology \emph{closed $p$-forms} seems to suggest,
\emph{being closed} is not a well defined property for a $p$-form, and
there is indeed an entire space of ``\emph{closing
structures}'' on a given $p$-form, namely the
homotopy fiber $K(w)$.  As a point in $K(w)$ consists of the data
needed to ``close'' the $p$-form $w$, we will call $K(w)$ the
\emph{space of keys of $w$}. For future reference we record this in
the following definition.

\begin{df}\label{d4}
For $A\in \cdga_{k}$, and $w\in \A^{p}(A,n)$, 
the space $K(w)$, is defined to be the homotopy fiber of the natural map
$\A^{p,cl}(A,n) \longrightarrow \A^{p}(A,n)$
taken at $w$. It is called the \emph{space of keys of $w$}. 
\end{df}

We have defined two $\s$-functors $\A^{p}(-,n), \quad \A^{p,cl}(-,n) :
\cdga_{k} \longrightarrow \mathbb{S}$ that we consider as derived
pre-stacks, $\dAff_{k}^{op} \longrightarrow \mathbb{S}$.  When $A \in
\cdga_{k}$ is viewed as a derived scheme $X=\Spec\, A$, we will
obviously write
$$
\A^{p}(X,n)=\A^{p}(A,n) \quad \A^{p,cl}(X,n)=\A^{p,cl}(A,n).
$$

\begin{prop}\label{p2}
The derived pre-stacks $\A^{p}(-,n)$ and $\A^{p,cl}(-,n)$ are derived
stacks for the \'etale topology.
\end{prop}
\textbf{Proof.} For $\A^{p}(-,n)$, we have by definition
$$
\A^{p}(\Spec\, A,n)\simeq
Map_{\dg_{k}}(k,(\wedge^{p}_{A}\mathbb{L}_{A/k})[n]).
$$ 
Therefore, the fact that $\A^{p}(-,n)$ satisfies \'etale descent
follows from the fact that the functor $\Spec\, A \mapsto
\wedge^{p}\mathbb{L}_{A/k}$ satisfies \'etale descent. This latter
$\s$-functor, when restricted to the small \'etale site of a derived
affine scheme $X$ is a quasi-coherent complex of
$\mathcal{O}_{X}$-modules, and thus is a derived stack for the \'etale
(and in fact the fpqc) topology (see \cite[Lemma 2.2.2.13]{hagII}).

In the same way, the derived pre-stack
$$ 
\Spec\, A \mapsto \mathbf{DR}(A/k)[n-p]\simeq
\bigoplus_{q}(\wedge^{q}_{A}\mathbb{L}_{A/k})[q-p+n] 
$$ 
is a quasi-coherent complex on the small \'etale site of $X=\Spec\,
A$, and so
satisfies \'etale descent. Therefore, Corollary \ref{cp1} implies that 
$$
\Spec\, A \mapsto NC^{w}(\mathbf{DR}(A/k)[n-p])\simeq
NC^{w}(\mathbf{DR}(A/k))[n-p] 
$$
is a derived stack for the \'etale topology. Taking the degree $p$
part and applying the 
Dold-Kan correspondence, we deduce that $\Spec\, A \mapsto
\A^{p,cl}(A,n)$ 
has descent for the \'etale topology. 
\hfill $\Box$

\

\smallskip

Proposition \ref{p2} enables us to globalize the definition of
$\A^{p}(X,n)$ and $\A^{p,cl}(X,n)$ to \emph{any} derived stack $X$, as
follows

\begin{df}\label{d5}
Let $F\in \dSt_{k}$ be a derived stack over $k$, $p$ and $n$ integers
with $p\geq 0$.  The \emph{space of $p$-forms, relative to $k$, of
degree $n$ on $F$} is defined by
$$
\A^{p}_{k}(F,n):=Map_{\dSt_{k}}(F,\A^{p}_{k}(-,n)).
$$
The \emph{space of closed $p$-forms, relative to $k$, of degree $n$ on
  $F$} is defined by 
$$
\A^{p,cl}_{k}(F,n):=Map_{\dSt_{k}}(F,\A^{p,cl}_{k}(-,n)).
$$
\end{df}

As before, when the base ring $k$ is clear from the context, we will not
include it in the notation.

Using the natural projection $\A^{p,cl}(-,n) \longrightarrow
\A^{p}(-,n)$, we have, for any $F$, a natural projection
$\A^{p,cl}(F,n) \longrightarrow \A^{p}(F,n)$. 

Definition \ref{d5} above has an alternative description, based on
Corollary \ref{cp1}. Consider the derived pre-stack $\Spec\, A \mapsto
\mathbf{DR}(A/k)$. It is a derived pre-stack with values in
$\medg_{k}$. As the forgetful $\s$-functor
$$\medg_{k} \longrightarrow \mdg_{k}$$ is conservative and preserves
$\s$-limits, Corollary \ref{cp1} implies that $\mathbf{DR}$ is a
derived stack with values in graded mixed complexes. By left Kan
extension, this derived stack extends uniquely to an $\s$-functor (see
e.g. \cite[\S 1.2]{chern} or \cite{lutop})
$$\mathbf{DR}(-/k) : \dSt_{k}^{op} \longrightarrow \medg_{k}.$$
Therefore, as in Definition 
\ref{d5}, we may give the following

\begin{df}\label{d5'}
For a derived stack $F$ we set
$$
NC^{w}(F/k):=NC^{w}(\mathbf{DR}(F/k)) \in \mdg_{k}.
$$
As usual, if $k$ is clear from the context, we will write
$$
\mathbf{DR}(F)=\mathbf{DR}(F/k) \qquad NC^{w}(F)=NC^{w}(F/k).
$$
\end{df}

Note that Corollary \ref{cp1} implies that we have natural equivalences
$$
\A^{p}_{k}(F,n)\simeq |\mathbf{DR}(F/k)[n-p](p)| \qquad
\A^{p,cl}_{k}(F,n)\simeq |NC^{w}(F/k)[n-p](p)|.
$$
The second of these equivalences is sometimes useful to compute
the spaces of closed $p$-forms by first computing explicitly the 
graded mixed complex $\mathbf{DR}(F/k)$ and then applying the
$\s$-functor $NC^{w}$. \\

The space of closed $p$-forms on a general derived stack $F$ can be a
rather complicated object, even when $F$ is a nice derived Artin
stack. The space of $p$-forms, however, has the following, expected,
description.

\begin{prop}\label{p3}
Let $F$ be a derived Artin stack over $k$, and $\mathbb{L}_{F/k} \in
L_{qcoh}(F)$ be its cotangent complex relative to $k$. Then, we have
an equivalence
$$
\A^{p}(F,n) \simeq
Map_{L_{qcoh}(F)}(\mathcal{O}_{F},\wedge^{p}\mathbb{L}_{F/k}[n]).
$$
This equivalence is functorial in $F$ in the obvious sense. In
particular, we have a functorial bijective map 
$$
\pi_{0}(\A^{p}(F,n)) \simeq H^{n}(F,\wedge^{p}\mathbb{L}_{F/k}).
$$
\end{prop}
\textbf{Proof.} We start by constructing a morphism
$$
\phi_{F} :
Map_{L_{qcoh}(F)}(\mathcal{O}_{F},\wedge^{p}\mathbb{L}_{F/k}[n])
\longrightarrow \A^{p}(F,n),
$$
functorial in $F$. By definition, the right hand side is given by 
$$
\A^{p}(F,n)\simeq \Holim_{X=\Spec\, A\in (\dAff/F)^{op}}\displaylimits
\A^{p}(X,n)\simeq \Holim_{X=\Spec\, A\in (\dAff/F)^{op}}\displaylimits
Map_{\dg_{k}}(k,\wedge^{p}_{A}\mathbb{L}_{A/k}[n]).
$$ 
On the $\s$-site $\dAff/F$, of derived affine schemes over $F$, we
have a (non-quasi-coherent) $\mathcal{O}$-module, denoted by
$\wedge^{p}\mathbb{L}$, and defined by
$$
X=\Spec\, A \mapsto \wedge^{p}_{A}\mathbb{L}_{A/k}.
$$
We also have the quasi-coherent $\mathcal{O}$-module
$\wedge^{p}\mathbb{L}_{F/k}$. There exists a natural  
morphism of $\mathcal{O}$-modules on $\dAff/F$, 
$\wedge^{p}\mathbb{L}_{F/k} \longrightarrow \wedge^{p}\mathbb{L}$.
obtained, over $u : \Spec\, A \longrightarrow F$, by the natural morphism
$u^{*}(\mathbb{L}_{F/k}) \longrightarrow \mathbb{L}_{A/k}$.
This defines a morphism on global sections
$$
Map_{L_{qcoh}(F)}(\mathcal{O}_{F},\wedge^{p}\mathbb{L}_{F/k}[n])
\longrightarrow  
\Holim_{X=\Spec\, A\in (\dAff/F)^{op}}
Map_{\dg_{k}}(k,\wedge^{p}_{A}\mathbb{L}_{A/k}[n]) \simeq \A^{p}(F,n).
$$
We will now check that this morphism
$$
\phi :
Map_{L_{qcoh}(F)}(\mathcal{O}_{F},\wedge^{p}\mathbb{L}_{F/k}[n])
\longrightarrow \A^{p}(F,n)
 $$ is an equivalence. For this we assume that $F$ is a derived Artin
stack which is $m$-geometric for some $m\geq 0$ (see \cite[\S
1.3.3]{hagII}), and we proceed by induction on $m$. For $m=0$, $F$ is
a derived affine stack and the statement is true. We assume that the
statement is correct for $(m-1)$-geometric derived stacks. We can
write $F$ as the quotient of a smooth Segal groupoid object $X_{*}$,
with $X_{i}$ being $(m-1)$-geometric for all $i$ (see \cite[\S
1.3.4]{hagII}).

We have a commutative square of descent maps
$$
\xymatrix{
Map_{L_{qcoh}(F)}(\mathcal{O}_{F},\wedge^{p}\mathbb{L}_{F/k}[n])
\ar[r] \ar[d]_-{\phi}&  
\Holim_{i\in
  \Delta}Map_{L_{qcoh}(X_{i})}(\mathcal{O}_{X_{i}},
\wedge^{p}(\mathbb{L}_{X_{i}/k})[n]) \ar[d]^-{\phi} \\  
\A^{p}(F,n)\ar[r] & \Holim_{i \in \Delta} \A^{p}(X_{i},n)}.
$$ 
By induction, the right vertical morphism between the homotopy
limits is an equivalence. The bottom horizontal morphism is also an
equivalence, because $F \mapsto \A^{p}(F,n)$ sends $\s$-colimits to
$\s$-limits by definition. It thus remains to show that the top
horizontal morphism is an equivalence. But this follows from the
following lemma. 

\begin{lem}\label{l1}
Let $X_{*}$ be a smooth Segal groupoid object in derived Artin stacks,
with quotient $F=|X_{*}|$. Then, for any integer $p$ and $n$, the
natural morphism
$$
Map_{L_{qcoh}(F)}(\mathcal{O}_{F},\wedge^{p}\mathbb{L}_{F/k}[n])
\longrightarrow 
\Holim_{i \in \Delta} \,
Map_{L_{qcoh}(X_{i})}(\mathcal{O}_{X_{i}},
\wedge^{p}\mathbb{L}_{X_{i}/k}[n])  
$$
is an equivalence. 
\end{lem}
\textit{Proof of lemma - } For any $n$, the morphisms appearing
in the statement of the lemma are retracts of the natural morphism
$$
Map_{L_{qcoh}(F)}\left(\mathcal{O}_{F},\bigoplus_{p}(
\wedge^{p}\mathbb{L}_{F/k})[n-p]\right)
\longrightarrow 
\Holim_{i \in
  \Delta}\left(Map_{L_{qcoh}(X_{i})}\left(\mathcal{O}_{X_{i}}, 
\bigoplus_{p}(\wedge^{p}\mathbb{L}_{X_{i}/k}) 
[n-p]\right)\right).
$$ 
The mapping spaces of the $\s$-category
$L_{qcoh}(F)$ are related to the derived Hom's by the formula
$$
Map_{L_{qcoh}(F)}(E,F) \simeq |\rh(E,F)|,
$$
and therefore it is enough  to prove that 
$$
\rh(\mathcal{O}_{F},\bigoplus_{p}(\wedge^{p}\mathbb{L}_{F/k})[-p])
\longrightarrow 
\Holim_{i \in
  \Delta}\rh(\mathcal{O}_{X_{i}},\bigoplus_{p}(\wedge^{p}
\mathbb{L}_{X_{i}/k})[-p])
$$
is a quasi-isomorphism of complexes of $k$-modules. 

For a derived stack $F$, we denote by $T^{1}(F)$ the shifted tangent
stack
$$
T^{1}(F):=\mathbf{Map}(\Spec\, k[e_{-1}],F),
$$ 
where $e_{-1}$ is in degree $-1$. The $\s$-functor $F \mapsto
T^{1}(F)$ preserves derived Artin stacks and finite homotopy
limits. Moreover, if $X \longrightarrow Y$ is a smooth and surjective
morphism of derived Artin stacks then $T^{1}(X) \longrightarrow
T^{1}(Y)$ is an epimorphism of derived stacks, as one may see by using
the infinitesimal criterion for smoothness \cite[\S 2.2.5]{hagII}. It
follows formally from these properties that the natural morphism
$|T^{1}X_*| \longrightarrow T^{1}|X_*|=T^{1}(F)$ is an
equivalence. This implies that we have
$$
\mathbb{R}\underline{Hom}_{L_{qcoh}(T^{1}(F))}(\mathcal{O}_{T^{1}(F)},
\mathcal{O}_{T^{1}(F)}) \simeq 
\Holim_{i\in \Delta} 
\mathbb{R}\underline{Hom}_{L_{qcoh}(T^{1}(X_i))}(\mathcal{O}_{T^{1}(X_i)},
\mathcal{O}_{T^{1}(X_i)}).
$$
Moreover, for any derived Artin stack $X$, $T^{1}(X)$ is affine over
$X$ and can be written  
as a relative spectrum (see \cite[Proposition 1.4.1.6]{hagII})
$$T^{1}(X)\simeq \Spec_{X}\, (Sym(\mathbb{L}_{X}[-1])).$$
In particular, we have natural quasi-isomorphisms
$$\rh(\OO_{T^{1}(X)},\OO_{T^{1}(X)}) \simeq
\rh(\OO_{X},\bigoplus_{p}\wedge^{p}_{\OO_{X}}(\mathbb{L}_{X})[-p])).$$

We thus deduce that the natural morphism
$$\mathbb{R}\underline{Hom}_{L_{qcoh}(F)}(\mathcal{O}_{F},
\bigoplus_{p}(\wedge^{p}\mathbb{L}_{F/k})[-p]) 
\longrightarrow 
\Holim_{i \in \Delta}\mathbb{R}
\underline{Hom}_{L_{qcoh}(X_{i})}(\mathcal{O}_{X_{i}},\bigoplus_{p} 
(\wedge^{p}\mathbb{L}_{X_{i}/k})[-p])$$
is a quasi-isomorphism, which implies the lemma. 
\hfill $\Box$ \\

This finishes the proof of the proposition.
\hfill $\Box$ \\

\begin{rmk}\label{rloop2}
The interpretations of $p$-forms and closed $p$-forms as functions and
invariant functions on the derived loop stacks given in \ref{rloop1}
has a global counterpart.  This globalization is not totally obvious,
as the construction $X \mapsto \mathcal{L}X$ is not compatible with
smooth gluing. However, it is possible to introduce the formal loop
stack $\mathcal{L}^{f}X$, the formal completion of $\mathcal{L}X$
along the constant loops $X \hookrightarrow \mathcal{L}X$, and to
prove that functions on $\mathcal{L}^{f}X$ have smooth descent (this
is essentially the same argument as in our Lemma \ref{l1}, see also
\cite{bena}). The group stack $\mathcal{H}$ acts on
$\mathcal{L}^{f}X$, and if we denote by
$$
q : [\mathcal{L}^{f}X/\mathcal{H}] \longrightarrow B\mathcal{H}
\longrightarrow B\mathbb{G}_{m}
$$
the natural projection, we have
$$
NC^{w}(F/k) \simeq q_{*}(\OO_{[\mathcal{L}^{f}X/\mathcal{H}]}).
$$
In other words, closed $p$-forms on $X$ can be interpreted 
as $B\mathbb{G}_{a}$-invariant (or equivalently 
$S^{1}$-invariant) functions on $\mathcal{L}^{f}X$.
\end{rmk}

\begin{rmk}\label{derhamd}{Let $F$ be a derived Artin stack locally of finite presentation
over $k$. If we consider the absolute $\epsilon$-differential $D_{E}$ (Remark \ref{epsilond}) for $E:=\mathbf{DR}(F/k)[n-p]=\mathbb{R}\Gamma(F, \oplus_{q\geq 0} (\wedge^{q}\mathbb{L}_{F})[q])[n-p]$, we get $$d_{DR}(p):=D_{E}(p): \mathrm{DR}(F/k)(p)[n-p]\longrightarrow \mathrm{NC}^{w}(F/k)(p+1)[n-(p+1)]$$ whose geometric realization gives the \emph{derived de Rham differential} (as a map of spaces) $$\textrm{d}_{DR}(p):=|d_{DR}(p)|: \mathcal{A}^{p}(F;n) \longrightarrow \mathcal{A}^{p+1,\, \textrm{cl}}(F;n).$$}
\end{rmk} 

\medskip

We are now ready to define the space of $n$\emph{-shifted symplectic
structures} on a derived Artin stack $F$.  To do this, let us recall
that if $F$ is a derived Artin stack (locally of finite presentation
over $k$, as are all of our derived Artin stacks by convention), then
the cotangent complex $\mathbb{L}_{F/k}$ is a perfect complex. In
particular it is a dualizable object in
$L_{qcoh}(F)$ (or equivalently dualizable in $D_{qcoh}(F)$, 
see \cite{chern}). Its dual will
be denoted by $\mathbb{T}_{F/k}$. Under this condition, any $2$-form
$\omega$ of degree $n$ on $F$, induces by Proposition \ref{p3} a
  morphism 
of quasi-coherent complexes
$$
\mathcal{O}_{F} \longrightarrow
(\wedge^{2}_{\mathcal{O}_{F}}\mathbb{L}_{F/k})[n],
$$
and thus, by duality, a morphism
$\wedge^{2}_{\mathcal{O}_{F}}\mathbb{T}_{F/k} \longrightarrow 
\mathcal{O}_{F}[n]$.
By adjunction, this induces a well defined morphism in $L_{qcoh}(F)$ 
$$
\Theta_{\omega} : \mathbb{T}_{F/k} \longrightarrow \mathbb{L}_{F/k}[n].
$$

\begin{df}\label{d6}
Let $F$ be a derived Artin stack, and $n\in \mathbb{Z}$.
\begin{enumerate}
\item 
A $2$-form $\omega \in \A^{2}_{k}(F,n)$ is \emph{non-degenerate} if 
the corresponding morphism in $D_{qcoh}(F)$
$$\Theta_{\omega} : \mathbb{T}_{F/k} \longrightarrow \mathbb{L}_{F/k}[n]$$
is an isomorphism. 

We denote by $\A^{2}_{k}(F,n)^{nd}$ the full subspace of
$\A^{2}_{k}(F,n)$ which is the union of all the connected components
consisting of non-degenerate $2$-forms of degree $n$ on $F$.

\item The \emph{space of $n$-shifted symplectic structures on $F$ (relative to
$k$)}, $\Symp_{k}(F,n)$ is defined by the following homotopy pull-back
square
$$
\xymatrix{
\Symp_{k}(F,n) \ar[r] \ar[d] & \A^{2,cl}_{k}(F,n) \ar[d] \\
\A^{2}_{k}(F,n)^{nd} \ar[r] & \A^{2}_{k}(F,n).}
$$
\end{enumerate}
We will simply write $\Symp(F,n)$ when the base ring $k$ is clear from
the context. 
\end{df}

Note that, by definition, $\Symp(F,n)$ is the full sub-space of
$\A^{2,cl}(F,n)$ defined by a unique condition on the underlying
$2$-form. \\

\begin{rmk}\label{r2}
As $F$ is locally of finite presentation, the cotangent complex
$\mathbb{L}_{F/k}$ is perfect and thus of bounded amplitude.  In
particular we see that at most one of the spaces $\Symp(F,n)$ can be
non-empty when $n$ varies in $\mathbb{Z}$, otherwise $\mathbb{L}_{F/k}$
will be periodic and thus not perfect.\\ More precisely, let us say that a derived Artin stack, locally of finite presentation over $k$,
has \emph{amplitude in} $[-m,n]$ (with $m, n\geq 0$) if its cotangent complex
has perfect amplitude in that range. Then such a derived stack might only carry
shifted symplectic
structures of degree $r=m-n$.
\end{rmk}

We conclude this section by describing the space of closed
$p$-forms in two simple situations: \emph{smooth schemes} (or more
generally Deligne-Mumford stacks), and \emph{classifying stacks} of
reductive group schemes. Finally we show that \emph{shifted cotangent stacks} carry a
canonical shifted symplectic structure. \\

\noindent \textbf{Smooth schemes.} We start with the case of a
smooth scheme $X$ over $\mathrm{Spec}\, k$. In this case
$\mathbb{L}_{X/k}\simeq \Omega_{X/k}^{1}$, and thus Proposition
\ref{p3} gives a description of the spaces of $p$-forms of degree $n$
as
$$
\A^{p}(X,n)\simeq |\mathbb{R}\Gamma(X,\Omega^{p}_{X/k})[n]|.
$$

Assume first that $X=\Spec\, A$ is smooth and affine. In this case we
know that $\mathbf{DR}(A/k)$ is naturally quasi-isomorphic, as a
graded mixed complex, to $DR(A/k)$, the usual, underived, de Rham
algebra of $A$ over $k$ (see \cite{tove}). We can then explicitly
compute the graded complex $NC^{w}(A/k)$. By applying directly the
definitions in Section \ref{gmc}, we get, for any $p\geq 0$
$$
NC^{w}(A/k)[-p](p)\simeq NC^{w}(DR(A/k))[-p](p) \simeq \Omega^{\geq
  p}_{A/k},
$$
where $\Omega^{\geq p}_{A/k}$ is the naively truncated de Rham complex
$$
\xymatrix{
\Omega^{p}_{A/k} \ar[r] & \Omega^{p+1}_{A/k} \ar[r] &
\Omega^{p+2}_{A/k} \ar[r] & \dots,}
$$ 
where $\Omega^{p}_{A/k}$ sits in degree $0$. From this we deduce
that the space of closed $p$-forms of degree $n$ on $\Spec\, A$ is
$$
\A^{p,cl}(\Spec\, A,n)\simeq |\Omega^{\geq p}_{A/k}[n]|.
$$
In particular, we have
$$
\pi_{i}(\A^{p,cl}(\Spec\, A,n)) \simeq \left\{ 
\begin{array}{ll} 
0 & if \; n<0 \\ H^{p+n-i}_{DR}(\Spec\, A) & if \; 0 \leq i < n \\
\Omega^{p,cl}_{A/k} & if \; i=n
\end{array} 
\right.
$$

By descent, we have similar formulas for a general smooth scheme $X$
(or more generally 
for a smooth Deligne-Mumford stack over $k$)
$$
\A^{p,cl}(X,n)\simeq |\mathbb{R}\Gamma(X,\Omega^{\geq p}_{X/k})[n]|,
$$
$$
\pi_{i}(\A^{p,cl}(X,n))\simeq H^{n-i}(X,\Omega^{\geq p}_{X/k}).
$$
As a consequence, we have the following three important properties for
a smooth scheme $X$ 
\begin{enumerate}
\item The space $\A^{p,cl}(X,0)$ of closed $p$-forms of degree $0$ is
equivalent to the discrete set $\Gamma(X,\Omega^{p,cl}_{X/k})$ of
\emph{usual} closed $p$-forms.
\item The spaces $\A^{p,cl}(X,n)$ are empty for $n<0$.
\item The spaces $\A^{p,cl}(X,n)$ are $n$-truncated for $n\geq 0$, and
we have
$$\pi_{i}(\A^{p,cl}(X,n))\simeq H^{n-i}(X,\Omega^{\geq p}_{X/k}).$$
\end{enumerate}

We can use these properties to describe $\Symp(X,n)$ for all smooth
schemes $X$ and all $n$. Indeed, property $(2)$ implies that
$\Symp(X,n)$ is empty for all $n<0$. Property $(1)$ gives that
$\Symp(X,0)$ is equivalent to the set of usual symplectic forms on
$X$.  Finally, $\Symp(X,n)$ is also empty for $n>0$, as all closed
$2$-forms of degree $n$ must be degenerate, since we cannot have
$\Omega_{X}^{1}[n]\simeq \mathrm{T}_{X}$ unless $n=0$ (note that both
$\Omega_{X}^{1}$ and $\mathrm{T}_{X}$ are complexes concentrated in
degree $0$ since $X$ is smooth).  In other words, on smooth schemes
there are no derived symplectic forms of degree $\neq 0$, and those in
degree $0$ are just the usual ones.\\

Note also an interesting consequence of the above properties. When $X$
is smooth and also \emph{proper} over $\Spec\, k$, then we have
$$
\pi_{i}(\A^{p,cl}(X,n))\simeq F^{p}H^{n+p-i}_{DR}(X/k),
$$ 
where $F^{*}$ stands for the Hodge filtration on the de Rham
cohomology $H^{*}_{DR}(X/k)$. \\

\noindent \textbf{Classifying stacks.} Let now $G$ be an \emph{affine
smooth group scheme} over $\mathrm{Spec}\, k$, and consider its
classifying stack $BG$, viewed as a derived Artin stack. If
$\mathfrak{g}$ denotes the Lie algebra of $G$, the cotangent complex
of $BG$ is $\mathfrak{g}^{\vee}[-1]$, considered as a quasi-coherent
complex on $BG$ via the adjoint action of $G$ on its lie algebra
$\mathfrak{g}$. In particular, using Proposition \ref{p3}, we find
that the spaces of $p$-forms on $BG$ are given by
$$ 
\A^{p}(BG,n)\simeq
|\,\mathbb{H}(G,Sym^{p}_{k}(\mathfrak{g}^{\vee}))[n-p]\,|,
$$ 
where $\mathbb{H}(G,-)$ denotes the Hochschild cohomology complex of the
affine group scheme $G$ with quasi-coherent coefficients (as in 
 \cite[\S 3.3]{dega} and \cite[\S 1.5]{to3}).  In particular, when
$G$ is \emph{reductive}, we find
$$
\A^{p}(BG,n)\simeq |\,Sym^{p}_{k}(\mathfrak{g}^{\vee})^{G}[n-p]\,|,
$$
and, equivalently, 
$$ 
\pi_{i}(\A^{p}(BG,n)) \simeq \left\{ \begin{array}{lc} 0 &
\textrm{if}\; i\neq n-p \\ Sym^{p}_{k}(\mathfrak{g}^{\vee})^{G} &
\textrm{if} \; i=n-p
\end{array}
\right.
$$ 
Still under the hypothesis that $G$ is reductive, let us now
compute the space $\A^{p,cl}(BG,n)$. For this we start by computing
the graded mixed complex $\mathbf{DR}(BG)$. Proposition \ref{p3} tells
us that the underlying graded complex of $\mathbf{DR}(BG)$ is
cohomologically concentrated in degree $0$
$$
\mathbf{DR}(BG) \simeq (Sym^{*}_{k}\mathfrak{g}^{\vee})^{G}[0].
$$ 
It follows that $\mathbf{DR}(BG)$, as a graded mixed complex is
quasi-isomorphic to $Sym^{*}_{k}(\mathfrak{g}^{\vee})[0]$, where the
$\epsilon$-action is trivial, and the grading is the natural grading on
$Sym^{*}_{k}(\mathfrak{g}^{\vee})$ (where $\mathfrak{g}^{\vee}$ is
assigned weight $1$). A direct consequence is that we have a natural
quasi-isomorphism of graded complexes
$$
NC^{w}(\mathbf{DR}(BG)) \simeq \bigoplus_{i \geq 0}
Sym^{*}_{k}(\mathfrak{g}^{\vee})^{G}[-2i].
$$
Thus, we have 
$$
\A^{p,cl}(BG,n)\simeq |\bigoplus_{i \geq
  0}Sym^{p+i}_{k}(\mathfrak{g}^{\vee})^{G}[n-p-2i]|.
$$
In particular
$$
\pi_{0}(\A^{p,cl}(BG,n))\simeq \left\{
\begin{array}{lc}
  0 & \textrm{if} \; n \; \textrm{is odd} \\
  Sym^{p}_{k}(\mathfrak{g}^{\vee})^{G} & \textrm{if} \; n \;
  \textrm{is even}
\end{array}
\right.  
$$

For degree reasons, a closed $2$-form of degree $n$ can be non-degenerate
on $BG$ only when $n=2$. Moreover, 
$$
\pi_{0}(\A^{2,cl}(BG,2))\simeq Sym^{2}_{k}(\mathfrak{g}^{\vee})^{G},
$$ 
is the $k$-module of $G$-invariant symmetric bilinear forms on
$\mathfrak{g}$. Such a closed $2$-form is non-degenerate if and only
if the corresponding bilinear form on $\mathfrak{g}^{\vee}$ is
non-degenerate in the usual sense. As a consequence, we have
$$ 
\pi_{0}(\Symp(BG,2)) \simeq \{\textrm{non-degenerate} \;
G\textrm{-invariant quadratic forms on} \; \mathfrak{g}\}.
$$ 
If $G$ is assumed to have simple
geometric fibers, then $\Symp(BG,2)$ possesses essentially a unique
element. Indeed, $Sym^{2}_{k}(\mathfrak{g}^{\vee})^{G}$ is a
projective module - being a direct factor in
$Sym^{2}_{k}(\mathfrak{g}^{\vee})$. Moreover, it is well known that
this projective module is of rank one when $k$ is a field, and this
implies, by base change, that $Sym^{2}_{k}(\mathfrak{g}^{\vee})^{G}$
is in fact a line bundle on $\Spec\, k$. Nowhere vanishing sections of
this line bundle corresponds to $2$-shifted symplectic forms on $BG$. We thus
see that, at least locally on the Zariski topology of $\Spec\, k$,
there is a $2$-shifted symplectic form on $BG$, which is unique up to a
multiplication by an invertible element of $k$.

When a reductive group scheme $G$ is realized as a closed subgroup
scheme of $\mathrm{GL}_{n}$, then there is a natural element in
$\Symp(BG,2)$. The inclusion $G \hookrightarrow \mathrm{GL}_n$ defines
a faithful representation $V$ of $G$ on $k^{n}$.  This representation
has a character, which is a $G$-invariant function on $G$. This
function can be restricted to the formal completion of $G$ at the
identity, to get a well defined element in
$$
\alpha_{V} \in \mathcal{O}(\widehat{G}^{G}) \simeq
\widehat{Sym}_{k}(\mathfrak{g}^{\vee})^{G}.
$$ 
The degree $2$ part of this element provides an invariant symmetric
bilinear form on $\mathfrak{g}$, which is non-degenerate because $G$
is reductive, and thus a $2$-shifted symplectic form on $BG$. In different
terms, $B\mathrm{GL}_{n}$ has a canonical $2$-shifted symplectic structure,
given by the bilinear form $(A,B) \mapsto Tr(AB)$, defined on the
$k$-modules $\mathcal{M}_{n}(k)$ of $n \times n$ matrices. The
inclusion $G \hookrightarrow \mathrm{GL}_{n}$ defines a morphism of
stacks $BG \longrightarrow B\mathrm{GL}_{n}$, and the pull-back of the
canonical $2$-shifted symplectic form on $BG$ remains a $2$-shifted symplectic form.\\

\noindent \textbf{Shifted cotangent stacks.} We define the $n$-shifted cotangent stack of a Deligne-Mumford stack, and prove that it carries a canonical $n$-shifted symplectic structure.

\begin{df} Let $n\in \mathbb{Z}$, and $X$ be a derived Artin stack locally of finite presentation over $k$. We define the \emph{$n$-shifted cotangent (derived) stack}  as $$\mathrm{T}^*X[n]:=\mathbb{R}\mathrm{Spec}\, \mathrm{Sym}_{\mathcal{O}_{X}}(\mathbb{T}_{X}[-n]).$$ For any $n$, we have a map of derived stacks $p[n]:\mathrm{T}^*X[n]\longrightarrow X$, induced by the canonical map $\mathcal{O}_{X} \longrightarrow \mathrm{Sym}_{\mathcal{O}_{X}}(\mathbb{T}_{X}[-n]).$
\end{df}

On $\mathrm{T}^*X[n]$ we have a \emph{canonical Liouville} $n$-shifted $1$-form. The idea - like in differential geometry, where there is obviously no shift - is that the pullback of $p[n]: T^{*}X[n] \rightarrow X$ along $p[n]$ itself has a canonical section, the diagonal, and any such section gives a horizontal $1$-form on $T^{*}X[n]$.\\
More precisely, we start by considering the inclusion $$\mathbb{T}_{X}[-n] \hookrightarrow \mathrm{Sym}_{\mathcal{O}_{X}}(\mathbb{T}_{X}[-n]).$$ Since $\mathbb{T}_{X}$ is perfect, this corresponds uniquely - by adjunction - to a map $$\mathcal{O}_{X} \rightarrow \mathbb{T}_{X}[-n]^{\vee} \otimes_{\mathcal{O}_{X}} \mathrm{Sym}_{\mathcal{O}_{X}}(\mathbb{T}_{X}[-n]) \simeq \mathbb{L}_{X}[n] \otimes_{\mathcal{O}_{X}} \mathrm{Sym}_{\mathcal{O}_{X}}(\mathbb{T}_{X}[-n])\simeq p[n]_{*}p[n]^{*}\mathbb{L}_{X}[n],$$ that, again by adjution, yields a map $$\mathcal{O}_{T^{*}X[n]}\rightarrow p[n]^{*}\mathbb{L}_{X}[n].$$ By composing this arrow with the shift-by-$n$ of the canonical map $p[n]^{*}\mathbb{L}_{X}\rightarrow \mathbb{L}_{T^{*}X[n]}$, we obtain the \emph{Liouville $n$-shifted $1$-form} on $\mathrm{T}^*X[n]$ $$\lambda(X;n) : \mathcal{O}_{T^{*}X[n]}\rightarrow \mathbb{L}_{T^{*}X[n]}[n].$$ Nothe that $\lambda$ is horizontal, by definition, i.e. the composite $\mathcal{O}_{T^{*}X[n]}\rightarrow \mathbb{L}_{T^{*}X[n]}[n] \rightarrow \mathbb{L}_{T^{*}X[n]/X}[n]$ is zero.\\


Recall from Remark \ref{derhamd}, the existence of a derived de Rham differential $$\textrm{d}_{DR}:=|d_{DR}(1)|: \mathcal{A}^{1}(T^{*}X[n];n) \longrightarrow \mathcal{A}^{2, \textrm{cl}}(T^{*}X[n];n).$$ Let's denote by $\textrm{d}_{DR}(\lambda(X;n)) \in \pi_{0}(\mathcal{A}^{2, \textrm{cl}}(T^{*}X[n];n)) = H^{n-2}(NC^{w}(T^{*}X[n])(2))$ the induced $n$-shifted closed $2$-form. 

\begin{prop}\label{yes} If $X$ is a derived Deligne-Mumford stack, then the underlying $n$-shifted $2$-form of $\omega:=\textrm{d}_{DR}(\lambda(X;n))$ is non degenerate, i.e. $\omega$ is symplectic.
\end{prop}

\noindent \textbf{Proof.} Let us simply denote by $d\lambda^{\flat}: \mathbb{T}_{\mathrm{T}^*X[n]} \rightarrow \mathbb{L}_{\mathrm{T}^*X[n]}[n]$ the map induced - by adjunction - by the $2$-form in $ H^{n}(T^{*}_{X}[n], \wedge^{2}\mathbb{L}_{\mathrm{T}^*X[n]})$ underlying $\textrm{d}_{DR}(\lambda(X;n))$. We want to prove that  $d\lambda^{\flat}$ is a quasi-isomorphism. Since \'etale maps induce  equivalences of cotangent complexes, by using an \'etale atlas $\{U_{i} = \mathbf{Spec}(A_{i}) \rightarrow X \}$ for $X$, and the induced \'etale atlas $\{\mathrm{T}^*U_{i}[n]= \mathbf{Spec}(B_{i}) \rightarrow \mathrm{T}^*X[n] \}$ on $\mathrm{T}^*X[n]$, it will suffice to prove the same statement upon restriction along any such \'etale map $\mathbf{Spec}(B_{i}) \rightarrow \mathrm{T}^*X[n]$. By naturality of the construction of $d\lambda^{\flat}$ with respect to \'etale maps, it will be enough to prove the proposition for $X= \mathbf{Spec}(A)$ (with $A$ quasi-free).


So, let $X=\mathbf{Spec}(A)$, where $A$ is a quasi-free $k$-cdga on the quasi-basis $\{ x_{i} \}_{i\in I}$. Then $$\mathbb{L}_{X}\simeq \Omega_{A/k}^{1} = \oplus_{i\in I} A\delta x_{i},$$  with $|\delta x_{i}|= |x_{i}|$, with the usual differential defined by $d(a\delta x_{i})= d(a)\delta x_{i} + (-1)^{|a|}a\delta (dx_{i})$, where $\delta$ is the unique derivation $A\rightarrow \Omega_{A/k}^{1}$ extending $x_{i} \mapsto \delta x_{i}$. Therefore, $$\mathbb{T}_{X} \simeq \mathrm{T}_{A/k} = \oplus_{i\in I} A \xi_{i},$$ where $\xi_{i}$ is dual to $\delta x_{i}$, and $|\xi_{i}|=-|x_{i}|$. Moreover $$Y:=T^*X[n]= \mathbf{Spec}(\mathrm{Sym_{A}}(\mathbb{T}_{X}[-n]))\simeq \mathbf{Spec} (B),$$ where $B$ is quasi-free over $A$ with quasi-basis $\{y_{i}\}_{i\in I}$, $|y_{i}|= - |x_{i}| +n$. 
In other words, $B$ is quasi-free over $k$ with quasi-basis $\{ x_{i}, y_{i}\}_{i\in I}$, $|y_{i}|= - |x_{i}| +n$. 
Therefore $$\mathbb{L}_{Y}\simeq \Omega_{B/k}^{1} = \oplus_{i\in I} B \delta x_{i} \oplus B\delta y_{i}, $$ with $|\delta x_{i}|= |x_{i}|$, $|\delta y_{i}|= -|x_{i}| +n$ (with its usual differential), and $$\mathbb{T}_{Y} \simeq \mathrm{T}_{B/k} = \oplus_{i\in I} B \xi_{i} \oplus B\eta_{i}, $$ where $\xi_{i}$  is dual to $\delta x_{i}$, $\eta_{i}$ is dual to $\delta y_{i}$, and $|\xi_{i}|= -|x_{i}|$, $|\eta_{i}|= |x_{i}| - n$. Moreover $$\mathbb{L}_{Y/X} \simeq \Omega^{1}_{B/A} = \oplus_{i \in I} B\delta y_{i}\, , \,\,\,\, \mathbb{T}_{Y/X} \simeq \mathrm{T}_{B/A} = \oplus_{i \in I} B\eta_{i}.$$ In these terms, the $n$-shifted Liouville $1$-form $\lambda$ on $Y$ is given by $\lambda = \sum_{i \in I}(-1)^{|y_i|} y_{i} \delta x_{i}$.
 Note that $\lambda$ is an element of degree $0$ in $DR(Y/k)(1)[n-1]\simeq \Omega_{B/k}^{1}[n]$. Now, by definition of the de Rham differential (Remark \ref{derhamd}), we have $$d_{DR}\lambda= (d_{DR}(\sum_{i \in I} (-1)^{|y_i|} y_{i} \delta x_{i})= \sum_{i \in I} (-1)^{|y_i|} \delta y_{i}\wedge \delta x_{i}, 0, \ldots, 0, \ldots) \in (NC^{w}(DR(Y/k))(2)[n-2])^{0},$$
so that the $n$-shifted $2$-form underlying $d_{DR}\lambda$ is $\sum_{i \in I} (-1)^{|y_i|} \delta y_{i}\wedge \delta x_{i}$.
In particular, we have 
$$d\lambda^{\flat}: \mathrm{T}_{B/k} =  \oplus_{i\in I} B \xi_{i} \oplus B\eta_{i} \longrightarrow \oplus_{i\in I} B \delta x_{i}[n] \oplus B\delta y_{i}[n] = \Omega^{1}_{B/k}[n],$$ $$\xi_{i} \longmapsto - (-1)^{|y_i|}\delta y_{i}[n] \, , \,\,\,\, \eta_{i} \longmapsto (-1)^{|y_i|} \delta x_{i}[n],$$ 
and this is, by inspection, an isomorphism of graded $B$-modules, and we conclude that $\omega$ is indeed symplectic (since we already knew that $d\lambda^{\flat}$ is a map of $B$-dg-modules).



\hfill $\Box$

\begin{rmk} It is very likely that Proposition \ref{yes} holds for any derived Artin stack, we just did not investigate how one may deal with smooth atlases instead of \'etale ones.
\end{rmk}

\section{Existence of shifted symplectic structures}

We prove in this section \emph{three} existence results for $n$-shifted symplectic forms
on certain derived stacks. We start by the \emph{mapping derived stack} 
of an oriented object to an $n$-shifted symplectic target, which is surely the most important 
of the existence theorems given below. We also introduce the notion of a \emph{Lagrangian 
structure} on a morphism with target an $n$-shifted symplectic derived stack, and 
show that the fibered product of two such morphisms comes equipped with a natural
$(n-1)$-shifted symplectic form. Finally, we explain a construction of 
$2$-shifted symplectic forms using the Chern character construction of \cite{chern}, and 
apply it to exhibit a natural $2$-shifted symplectic structure on $\Parf$, the
derived moduli stack of perfect complexes. \\

\subsection{Mapping stacks}

The following construction is well-known in differential geometry and provides, together with \cite{aksz}, heuristics for our constructions below. Let $M$ be a compact $C^{\infty}$-manifold of dimension $m$, $N$ be a $C^{\infty}$-manifold, and $Map_{C^{\infty}}(M,N)$ the Fr\'echet manifold of $C^{\infty}$-maps from $M$ to $N$. The following diagram $$\xymatrix{ & M\times Map_{C^{\infty}}(M,N) \ar[ld]^-{ev} \ar[rd]_-{pr_{M}} & \\ N & & M}$$
induces a map $$\Omega_{M}^{p} \times \Omega_{N}^{q} \rightarrow \Omega_{Map_{C^{\infty}}(M,N)}^{p+q-m} : (\alpha, \beta) \mapsto \int_{M} pr^{*}_{M}\alpha \wedge ev^{*} \beta := \widehat{\alpha \beta},$$ where $\int_{M}$ denotes integration along the fiber, sometimes called the \emph{hat-product} (see e.g. \cite{viz}). Now, if $(N,\omega)$ is symplectic, then $\widehat{1\omega} \in \Omega_{Map(M,N)}^{2-m}$ defines a symplectic form on $Map_{C^{\infty}}(M,N).$\\ 
In the derived algebraic geometry setting we are concerned with in this paper, we will need a replacement for Poincar\'e duality and for the notion of an orientation. The first one is given by Serre duality (in the more general context of Calabi-Yau categories), while the second one will be that of $\OO$\emph{-orientation} (Definition \ref{d9}). The notion of $\OO$-orientation will allow for a quasi-coherent variant of \emph{integration along the fiber} (Definition \ref{d8}).\\

\noindent We start by considering the following finiteness conditions on 
derived stacks.

\begin{df}\label{d7}
A derived stack $X$, over a derived affine scheme  $\Spec\, A$, 
is \emph{strictly $\OO$-compact over $A$} if it satisfies the following two 
conditions
\begin{enumerate}
\item $\OO_{X}$ is a compact object in $D_{qcoh}(X)$.
\item For any perfect complex $E$ on $X$, the $A$-dg-module
$$C(X,E) := \rh(\OO_{X},E)$$
is perfect. 
\end{enumerate}
A derived stack $X$ over $k$ is \emph{$\OO$-compact} if for any derived affine scheme $\Spec\, A$
the derived stack $X\times \Spec\, A$ is strictly $\OO$-compact over $A$.
\end{df}

\begin{rmk}\label{r3}
Since perfect complexes are exactly the dualizable objects in $D_{qcoh}(X)$, condition 
$(1)$ of the definition above implies that all perfect complexes are 
compact in $D_{qcoh}(X)$ for an $\OO$-compact derived stack $X$, as well as 
all perfect complexes on $X\times \Spec\, A$ for any $A \in \ncdga_{k}$. 
\end{rmk}

The main property of (strictly) $\OO$-compact 
derived stacks is the existence, for any other derived stack $F$, of
a morphism of graded mixed complexes (over $k$)
$$\kappa_{F,X} : \mathbf{DR}(F\times X) \longrightarrow \mathbf{DR}(F)\otimes_{k} C(X,\OO_{X}),$$
functorial in $F$, where $C(X,\OO_{X})$ is considered pure of weight $0$ with trivial
mixed structure. It is defined as follows. 

Since $X$ is, by hypothesis, $\OO$-compact, the complex $C(X,\OO_{X})$ is perfect over $k$, so 
the $\s$-endofunctor
$E \mapsto E\otimes_{k} C(X,\OO_{X})$, of the $\s$-category of mixed graded complexes, commutes with $\s$-limits. 
Also, by definition, the functor $F \mapsto \mathbf{DR}(F\times X)$ sends $\s$-colimits to $\s$-limits. Therefore, 
the two $\s$-functors
$$\mathbf{DR}(-\times X), \,\,\,  \mathbf{DR}(-)\otimes_{k} C(X,\OO_{X}) : \dSt_{k}^{op} \longrightarrow 
\medg_{k}$$
send $\s$-colimits of derived stacks to $\s$-limits of graded complexes. By left Kan extensions (see 
\cite{lutop} or \cite[\S 1.2]{chern}), 
we see that in order to construct a natural transformation
$\kappa_{-,X} : \mathbf{DR}(-\times X) \longrightarrow \mathbf{DR}(-)\otimes_{k} C(X,\OO_{X})$,
it is enough to construct a natural transformation between the two $\s$-functors restricted
to the $\s$-category of derived affine schemes. By definition, these two $\s$-functors, restricted
to derived affine schemes are given as follows
$$\begin{array}{cccc}
\mathbf{DR}(-) \otimes_{k} C(X,\OO_{X}) : & \dAff_{k}^{op}\equiv\ncdga_{k} & \longrightarrow & \medg_{k} \\
 & A & \longmapsto & \mathbf{DR}(A)\otimes_{k} C(X,\OO_{X})
 \end{array}$$
$$\begin{array}{cccc}
\mathbf{DR}(-\times X) : & \dAff_{k}^{op}\equiv\ncdga_{k} & \longrightarrow & \medg_{k} \\
 & A & \longmapsto & H(X,\mathbf{DR}(A\otimes_{k}\OO_{X})),
 \end{array}$$
where $\mathbf{DR}(A\otimes_{k}\OO_{X})$ is the stack of mixed graded complexes on $X$
sending $\Spec B \rightarrow X$ to $\mathbf{DR}(A\otimes_{k}B)$, and $H(X,\mathbf{DR}(A\otimes_{k}\OO_{X}))$
denotes its global sections. 
 
For any two objects $B,C \in \ncdga$, we have a natural equivalence of graded mixed complexes
(Kunneth formula) 
$$\mathbf{DR}(B)\otimes_{k} \mathbf{DR}(C) \simeq \mathbf{DR}(B\otimes_{k}C),$$
induced by the identification 
$$\mathbb{L}_{B\otimes_{k} C} \simeq (\mathbb{L}_{B}\otimes_{k}C) \oplus (B\otimes_{k}\mathbb{L}_{C}).$$
Therefore, the $\s$-functor $\mathbf{DR}(-\times X)$ sends $A$ to 
$H(X,\mathbf{DR}(A)\otimes_{k}\mathbf{DR}(\OO_{X}))$. We
consider the natural projection on the component of weight zero (with 
trivial mixed structure) $\mathbf{DR}(\OO_{X}) \longrightarrow \OO_{X}$, and obtain a morphism
$\mathbf{DR}(-\times X) \longrightarrow H(X,\mathbf{DR}(-)\otimes_{k} \OO_{X})\simeq
C(X,\mathbf{DR}(-)\otimes_{k} \OO_{X})$. As
$X$ is $\OO$-compact, $C(X,-)$ commutes with colimits of quasi-coherent sheaves, and thus the natural morphism
$$\mathbf{DR}(-)\otimes_{k} C(X,\OO_{X}) \longrightarrow C(X,\mathbf{DR}(-)\otimes_{k} \OO_{X})$$
is an equivalence. 

We thus have defined a natural transformation of $\s$-functors
$$\mathbf{DR}(-\times X) \longrightarrow C(X,\mathbf{DR}(-)\otimes_{k} \OO_{X}) \simeq \mathbf{DR}(-) 
\otimes_{k}C(X,\OO_{X}),$$
which defines our morphism of graded mixted complexes
$$\kappa_{F,X} : \mathbf{DR}(F\times X) \longrightarrow \mathbf{DR}(F)\otimes_{k} C(X,\OO_{X}).$$
  
We can apply the $\s$-functor $NC^w$ to the morphism above, in order to get  
a morphism of graded complexes
$$\kappa_{F,X} : NC^{w}(F\times X) \longrightarrow NC^{w}(\mathbf{DR}(F)\otimes_{k} C(X,\OO_{X})).$$
As $C(X,\OO_{X})$ is a perfect complex over $k$, the morphism
$$NC^{w}(\mathbf{DR}(F))\otimes_{k} C(X,\OO_{X}) \longrightarrow NC^{w}(\mathbf{DR}(F)\otimes_{k} C(X,\OO_{X}))$$
is an equivalence of graded complexes. 

To summarized, we have defined for any derived stack $F$ and any $\OO$-compact derived stack $X$, 
a commutative square of graded complexes

$$\xymatrix{
NC^{w}(F\times X) \ar[r]^-{\kappa_{F,X}} \ar[d] & NC^{w}(F) \otimes_{k}C(X,\OO_{X}) \ar[d] \\
\mathbf{DR}(F\times X) \ar[r]_-{\kappa_{F,X}} & \mathbf{DR}(F) \otimes_{k} C(X,\OO_{X}),}$$
where the vertical morphisms are the projections $NC^{w} \longrightarrow \mathbf{DR}$. 

We keep the hypothesis that $X$ is an $\OO$-compact derived stack, and we assume further that, for some 
integer $d\in \mathbb{Z}$, 
we are given a map 
$$\eta : C(X,\OO_{X}) \longrightarrow k[d],$$
in the derived category $D(k)$. Then, for any derived stack 
$F$ we have a natural morphism of graded complexes
$$\xymatrix{
NC^{w}(F\times X) \ar[r]^-{\kappa_{F,X}} & NC^{w}(F)\otimes_{k} C(X,\OO_{X})
 \ar[r]^-{id\otimes \eta} & NC^{w}(F)[d].}$$
This morphism, well defined in the homotopy category of graded complexes $\ho(\mdg_{k})$,
is called the \emph{integration along $\eta$}. 
 
\begin{df}\label{d8}
Let $F$ and $X$ be derived stacks, with $X$ $\OO$-compact, and 
let \linebreak $\eta : C(X,\OO_{X}) \longrightarrow k[d]$ be a morphism in $D(k)$,
for some integer $d$.
The \emph{integration map along $\eta$} is the morphism 
$$\int_{\eta} : \xymatrix{
NC^{w}(F\times X) \ar[r]^-{\kappa_{F,X}} & NC^{w}(F)\otimes_{k} C(X,\OO_{X})
\ar[r]^-{id\otimes \eta} & NC^{w}(F)[d]}$$ constructed above.
\end{df}

We also have a similar
morphism on the level of de Rham complexes
$$\int_{\eta} : \xymatrix{
\mathbf{DR}(F\times X) \ar[r]^-{\kappa_{F,X}} & \mathbf{DR}\otimes_{k} C(X,\OO_{X})
 \ar[r]^-{id\otimes \eta} & \mathbf{DR}(F)[d],}$$
in a way that we have a commutative diagram of mixed graded complexes 
$$\xymatrix{
NC^{w}(F\times X) \ar[d]\ar[r]^-{\int_{\eta}}  & NC^{w}(F)[d] \ar[d] \\
\mathbf{DR}(F\times X) \ar[r]^-{\int_{\eta}}  & \mathbf{DR}(F)[d].}$$

Let $X$ be an $\OO$-compact derived stack. The complex
$C(X,\OO_{X})$ possesses a natural structure of a (commutative) dg-algebra over $k$, 
so comes equipped with a cup-product morphism
$$C(X,\OO_{X}) \otimes_{k} C(X,\OO_{X}) \longrightarrow C(X,\OO_{X}).$$
In particular, any morphism $\eta : C(X,\OO_{X}) \longrightarrow k[-d]$ provides a morphism
$$- \cap \eta : \xymatrix{
C(X,\OO_{X}) \otimes C(X,\OO_{X}) \ar[r]^-{\cap} & C(X,\OO_{X}) 
\ar[r]^-{\eta} & k[-d].}$$
If we denote by $C(X,\OO_{X})^{\vee}:=\rh (C(X,\OO_{X}),k)$ the derived
dual of $C(X,\OO_{X})$, the morphism above defines an adjoint morphism
$$ - \cap \eta : C(X,\OO_{X}) \longrightarrow C(X,\OO_{X})^{\vee}[-d].$$
More generally, if $E$ is a perfect complex on $X$, of dual 
$E^{\vee}:=\rch (E,\OO_{X})$, the natural pairing
$$C(X,E) \otimes_{k} C(X,E^{\vee}) \longrightarrow C(X,\OO_{X}),$$
composed with $\eta$ induces a morphism
$$- \cap \eta  : C(X,E) \longrightarrow C(X,E^{\vee})^{\vee}[-d].$$

If moreover $A \in \ncdga_{k}$, the same is true for the derived
$A$-scheme $X_{A}:=X\times \Spec\, A$. The morphism $\eta$ induces a morphism
$$\eta_{A}=\eta \otimes id_{A} : C(X_{A},\OO_{X_{A}}) \simeq C(X,\OO_{X})\otimes_{k}A \longrightarrow
k[-d]\otimes_{k}A,$$
and for any perfect complex $E$ on $X_{A}$, 
the morphism $\eta_{A}$ induces a natural morphism
$$-\cap \eta_{A} : C(X_{A},E) \longrightarrow C(X_{A},E^{\vee})^{\vee}[-d],$$
where now $C(X_{A},E^{\vee})^{\vee}$ is the derived $A$-dual 
of $C(X_{A},E^{\vee})$. 

\begin{df}\label{d9}
Let $X$ be an $\OO$-compact derived stack and $d\in \mathbb{Z}$.  An \emph{$\OO$-orientation of degree $d$ on 
$X$} consists of a morphism of complexes
$$[X] : C(X,\OO_{X}) \longrightarrow k[-d],$$ 
such that for any $A\in \ncdga_{k}$ and any perfect complex $E$ on $X_{A}:=X\times \Spec\, A$, 
the morphism
$$- \cap [X]_{A} : C(X_{A},E) \longrightarrow C(X_{A},E^{\vee})^{\vee}[-d]$$
is a quasi-isomorphism of $A$-dg-modules.
\end{df}

We are now ready to state and prove our main existence statement. 

\begin{thm}\label{t1}
Let $F$ be a derived Artin stack  equipped with
an $n$-shifted symplectic form $\omega \in \Symp(F,n)$. Let 
$X$ be an $\OO$-compact derived stack equipped with an $\OO$-orientation $[X] : C(X,\OO_{X}) \longrightarrow k[-d]$ of degree
$d$. Assume that the derived mapping stack $\mathbf{Map}(X,F)$ is itself a derived Artin stack 
locally of finite presentation over $k$. Then, 
$\mathbf{Map}(X,F)$ carries a canonical $(n-d)$-shifted symplectic structure.  
\end{thm}

\textbf{Proof.} We let 
$\pi : X\times \mathbf{Map}(X,F) \longrightarrow F$ 
be the evaluation morphism. We have 
$\omega \in \Symp(F,n) \subset \A^{2,cl}(F,n)$, and this corresponds to 
a morphism of graded complexes
$$\omega : k[2-n](2) \longrightarrow NC^{w}(F).$$
Using the integration along the orientation $[X]$ of definition \ref{d8}, we consider the composition
$$\int_{[X]}\omega : \xymatrix{ k[2-n](2) \ar[r]^-{\omega} & NC^{w}(F) \ar[r]^-{\pi^{*}}
& NC^{w}(X\times \mathbf{Map}(X,F)) \ar[r]^-{\int_{[X]}} & NC^{w}(\mathbf{Map}(X,F))[-d].}$$
This is, by definition, a closed $2$-form of degree $(n-d)$ on $\mathbf{Map}(X,F)$. i.e.
$$\int_{[X]}\omega \in \A^{2,cl}(\mathbf{Map}(X,F),n-d).$$

It remains to show that this $2$-form is non-degenerate. For this, we have to determine the
underlying $2$-form of degree $(n-d)$. It is given by the following morphism 
$$\xymatrix{ k[2-n](2) \ar[r]^-{\omega} & \mathbf{DR}(F) \ar[r]^-{\pi^{*}}
& \mathbf{DR}(X\times \mathbf{Map}(X,F)) \ar[r]^-{\int_{[X]}} & \mathbf{DR}(\mathbf{Map}(X,F))[-d].}$$
If we unravel the definition of $\int_{[X]}$, we see that this $2$-form can be described as follows. 

First of all, let $x : \Spec\, A \longrightarrow \mathbf{Map}(X,F)$ be an $A$-point corresponding to a morphism of derived 
stacks
$$f : X\times \Spec\, A \longrightarrow F.$$
If $\mathbb{T}_{F}$ denotes the tangent complex of $F$, the tangent complex of $\mathbf{Map}(X,F)$ at the point $x$ is given by 
$$\mathbb{T}_{x}\mathbf{Map}(X,F) \simeq C(X_{A},f^{*}(\mathbb{T}_{F})).$$
The $2$-form $\omega$ defines a non-degenerate pairing of perfect complexes on $F$ 
$$\mathbb{T}_{F} \wedge \mathbb{T}_{F} \longrightarrow \OO_{F}[n],$$ 
which induces an alternate pairing of $A$-dg-modules
$$C(X_{A},f^{*}(\mathbb{T}_{F})) \wedge 
C(X_{A},f^{*}(\mathbb{T}_{F})) \longrightarrow C(X_{A},\OO_{X_{A}}[n]).$$
We can 
compose with the orientation $[X_{A}]$ to get a pairing of perfect $A$-dg-modules. 
$$C(X\times \Spec\, A,f^{*}(\mathbb{T}_{F})) \wedge 
C(X\times \Spec\, A,f^{*}(\mathbb{T}_{F})) \longrightarrow A[n-d].$$
By inspection, this pairing is the one induced by the $2$-form underlying
$\int_{[X]}\eta$. The fact that it is non-degenerate then follows from 
the definition of an orientation and the fact that $\omega$ is non-degenerate.
\hfill $\Box$ \\

\noindent Here follow some examples of derived stacks $X$ satisfying the condition of the theorem \ref{t1}.

\begin{itemize}
\item \textbf{Mapping stacks with Betti source}. Let $M$ be a compact, connected, and oriented topological manifold. We consider
$X=S(M)$ its simplicial set of singular simplices, as a constant derived stack. 
The category $D_{qcoh}(X)$ is then naturally equivalent to 
$D_{loc}(M,k)$, the full sub-category of the derived category of sheaves of
$k$-modules on the space $M$, consisting of objects with locally constant 
cohomology sheaves. In particular, we have functorial isomorphisms
$$H^{*}(X,E) \simeq H^{*}(M,\mathcal{E}),$$
for any $E \in D_{qcoh}(X)$ whose corresponding 
complex of $k$-modules on $M$ is denoted by $\mathcal{E}$. 
Perfect complexes on $X$ correspond to objects in $D_{loc}(M,k)$
locally quasi-isomorphic to bounded complexes of constant sheaves of 
projective modules of finite type.
This implies that $X$ is $\OO$-compact. Moreover, 
the orientation on $M$
determines a well defined fundamental class $[M] \in H_{d}(M,k)$, and thus a morphism
$[X] : C(M,k) \longrightarrow k[-d]$, where $d=\dim \, M$. Poincar\'e duality on $M$ 
implies that $[X]$ is an $\OO$-orientation on $X$. Finally, $M$ has the homotopy
type of a finite CW complex, so $X$ is a finite homotopy type. This implies that 
for any derived Artin stack $F$, $\mathbf{Map}(X,F)$ is a finite homotopy limit of copies of 
$F$, and thus is itself a derived Artin stack.

\item \textbf{Mapping stacks with de Rham source}. Let $Y$ be a smooth and proper Deligne-Mumford stack over $\mathrm{Spec}\, k$, 
with connected geometric fibers.  
Recall from \cite[Corollary 2.2.6.15]{hagII} that we can define a derived stack 
$Y_{DR}$, such that 
$$Y_{DR}(A):=Y(\pi_{0}(A)_{red}),$$
for any $A\in \ncdga_{k}$. We set $X:=Y_{DR}$. We have a natural equivalence
between $D_{qcoh}(X)$ and the derived category of $\mathcal{D}_{Y/k}$-modules with
quasi-coherent cohomology. Moreover,  
perfect complexes on $X$ correspond to complexes of $\mathcal{D}_{Y/k}$-modules whose
underlying quasi-coherent complexes are perfect over $Y$. In particular, if $E$ is
a perfect complex on $X$, corresponding to a complex of $\mathcal{D}_{Y/k}$-modules $\mathcal{E}$
perfect over $Y$, then we have
$$H^{*}(X,E)\simeq H_{DR}^{*}(Y/k,\mathcal{E}).$$
It follows easily that $X$ is $\OO$-compact. Moreover, the choice of 
a fundamental class in de Rham cohomology $[Y] \in H^{2d}_{DR}(Y/k,\mathcal{O})$
(where $d$ is the relative dimension of $Y$ over $\mathrm{Spec}\, k$)
determines a morphism
$$[X] : C(X,\OO) \longrightarrow k[-2d]$$
which, by Poincar\'e duality in de Rham cohomology, is an $\OO$-orientation on $X$.

Finally, the fact that $\mathbf{Map}(X,F)$ is a derived Artin stack when $F$ is one, 
can be deduced from Lurie's version of Artin representability criterion. We will be mainly 
interested in the special case where $F$ is either a \emph{smooth quasi-projective variety}, or
a \emph{classifying stack} $BG$, or the \emph{derived stack of perfect complexes} $\Parf$. In all these
specific situations, the fact that $\mathbf{Map}(X,F)$ is a derived Artin stack
locally of finite presentation can be found in \cite[\S 2.2.6.3]{hagII}
and \cite{si2}. 

\item \textbf{Mapping stacks with Dolbeault source}. The previous example has the following Dolbeault analog. Let 
again $Y$ be a smooth and proper Deligne-Mumford stack over $\Spec\, k$, 
with connected geometric fibers. We define $Y_{Dol}$ by (see \cite[\S 2]{si2})
$$Y_{Dol}:=B\widehat{\mathbb{T}}_{Y/k} \longrightarrow Y,$$
the classifying stack of the formal tangent bundle of $Y$ relative to $k$. 
We set $X:=Y_{Dol}$.
We know that a quasi-coherent complex $E$ on $X$ consists of a pair $(\mathcal{E}, \phi)$ where $\mathcal{E}$
is a quasi-coherent complex  on $Y$, and $\phi$ is a Higgs field $\phi$ on $\mathcal{E}$ (i.e. an 
action of the $\OO_{Y}$-algebra $Sym_{\OO_{Y}}(\mathbb{T}_{Y/k})$). Under this correspondence, we have
$$H^{*}(X,E) \simeq H^{*}_{Dol}(Y,\mathcal{E}).$$
This implies that $X$ is $\OO$-compact. As above, the choice of 
a fundamental class in Hodge cohomology
$[Y] \in H^{2d}_{Dol}(Y,\OO)\simeq H^{d}(Y,\Omega^{d}_{Y/k})$ determines a morphism
$$[X] : C(X,\OO) \longrightarrow k[-2d]$$
which, by Poincar\'e duality in Dolbeault cohomology, is indeed an $\OO$-orientation on $X$.

Again, the fact that $\mathbf{Map}(X,F)$ is a derived Artin stack when $F$ is one, 
can be deduced from Lurie's version of Artin representability criterion. We will be mainly 
interested in the special case where $F$ is either a \emph{smooth quasi-projective variety}, or
a \emph{classifying stack} $BG$, or the \emph{derived stack of perfect complexes} $\Parf$. In all these
specific situations, the fact that $\mathbf{Map}(X,F)$ is a derived Artin stack
locally of finite presentation can be found in \cite[\S 2.2.6.3]{hagII}
and \cite{si2}.

The Dolbeault and de Rham complexes can also be considered together at the same time, by taking $X
:= Y_{Hod} \to \mathbb{A}^{1}$, see \cite{si1}. More generally, we could
take $X$ to be any nice enough formal groupoid (as in \cite{si1}).

\item \textbf{Mapping stacks with Calabi-Yau source}. Let now $X$ be a smooth and proper Deligne-Mumford stack over $\mathrm{Spec}\, k$ of relative
dimension $d$, with connected geometric fibers. We assume that we are given an isomorphism of line bundles
$$u : \omega_{X/k}=\wedge^{d}\Omega^{1}_{X/k} \simeq \OO_X \equiv \OO.$$
Considered as a derived stack $X$, is automatically $\OO$-compact. Moreover, the isomorphism 
$u$, together with the trace map, defines an isomorphism
$$\xymatrix{H^{d}(X,\OO) \ar[r]^-{u} & H^{d}(X,\omega_{X/k}) \ar[r]^-{tr} & k.}$$
This isomorphism induces a well defined morphism of complexes
$$C(X,\OO) \longrightarrow k[-d],$$
which, by Serre duality, is indeed an $\OO$-orientation on $X$.

As above, the fact that $\mathbf{Map}(X,F)$ is a derived Artin stack when $F$ is one, 
can be deduced from Lurie's version of Artin representability criterion. We will be mainly 
interested in the special case where $F$ is either a \emph{smooth quasi-projective variety}, or
a \emph{classifying stack} $BG$, or the \emph{derived stack of perfect complexes} $\Parf$. In all these
specific situations, the fact that $\mathbf{Map}(X,F)$ is a derived Artin stack
locally of finite presentation can be found in \cite[\S 2.2.6.3]{hagII}
and \cite{tv}.

\end{itemize}

\bigskip

We gather the following consequences of Theorem \ref{t1} and of the examples above, in the following 

\begin{cor}\label{ct1}
Let $G$ be a reductive affine group scheme over $\mathrm{Spec}\, k$. Let $Y$ be a smooth and proper
Deligne-Mumford stack over $\mathrm{Spec}\, k$ with connected geometric fibers of relative dimension $d$.
Assume that we have fixed a non-degenerate $G$-invariant symmetric bilinear form on $\mathfrak{g}$. 
\begin{enumerate}
\item The choice of a fundamental class $[Y] \in H^{2d}_{DR}(Y,\OO)$
determines a canonical $2(1-d)$-shifted symplectic form on the derived stack 
 $$\mathbb{R}\textrm{Loc}_{DR}(Y,G):=\mathbf{Map}(Y_{DR},BG)$$
 of flat $G$-bundles on $Y$.
\item The choice of a fundamental class $[Y] \in H^{2d}_{Dol}(Y,\OO)$
determines a canonical $2(1-d)$-shifted symplectic form on the derived stack  
$$\mathbb{R}\textrm{Loc}_{Dol}(Y,G):=\mathbf{Map}(Y_{Dol},BG)$$
of Higgs $G$-bundles on $Y$.
\item When it exists, the choice of a trivialization isomorphism 
$\omega_{Y/k} \simeq \OO_{Y}$,  
determines a canonical $(2-d)$-shifted symplectic form on the derived stack 
of $G$-bundles on $Y$ 
$$\mathbb{R}\textrm{Bun}(Y,G):=\mathbf{Map}(Y,BG)$$
of $G$-bundles on $Y$.
\item If $M$ is a compact, orientable topological manifold of dimension $d$, then a choice
of a fundamental class $[M]\in H_{d}(M,k)$ determines a canonical 
$(2-d)$-shifted symplectic form on the derived stack 
$$\mathbb{R}\textrm{Loc}(M,G):=\mathbf{Map}(M,BG)$$
of local systems of principal $G$-bundles on $M$.
\end{enumerate}
\end{cor}

In section 3 we will explain how these $n$-shifted symplectic structures 
compare with the well known symplectic structures on certain coarse moduli spaces
(e.g. on character varieties, moduli spaces of stable sheaves on K3-surfaces, etc.).

\subsection{Lagrangian intersections}

We will be interested here in the study of derived symplectic structures induced on fiber products
of derived Artin stacks. In order to do this, 
we first need to introduce the notion of \emph{isotropic} and \emph{Lagrangian data} on a morphism with 
symplectic target, extending the usual notions
of isotropic and Lagrangian sub-varieties of a smooth symplectic manifold. We will 
show (Theorem \ref{t2}) that the fiber product of two morphisms with Lagrangian structures towards 
an $n$-shifted symplectic target is naturally equipped with an $(n-1)$-shifted symplectic structure. In particular, the derived 
intersection of two smooth (algebraic) usual Lagrangians in a smooth (algebraic) symplectic manifold carries 
a canonically induced $(-1)$-shifted symplectic structure. \\

We fix a derived Artin stack $F$ 
and an $n$-shifted symplectic form $\omega \in \Symp(F,n)$ on $F$. For $f : X \longrightarrow F$ a morphism of 
derived Artin stack, we have the pull-back closed $2$-form
$$f^{*}(\omega) \in \A^{2,cl}(X,n).$$

\begin{df}\label{d10}
An \emph{isotropic structure on $f$ (relative to $\omega$)} is a path 
between $0$ and $f^{*}(\omega)$ in the space $\A^{2,cl}(X,n)$. The space
of \emph{isotropic structures on $f$ (relative to $\omega$)} is defined to be the path space 
$$\mathsf{Isot}(f,\omega):=Path_{0,f^{*}(\omega)}(\A^{2,cl}(X,n)).$$
\end{df}

A Lagrangian structure on the morphism $f$ (with respect to $\omega$) will be
an isotropic structure satisfying some non-degeneracy condition. To introduce  
this condition, let us fix an isotropic structure $h \in \mathsf{Isot}(f,\omega)$. 
We consider the $2$-form $\mathbb{T}_{F} \wedge \mathbb{T}_{F} \longrightarrow \OO_{F}[n]$ underlying $\omega$,
as well as its pull-back on $X$, 
$f^{*}(\mathbb{T}_{F}) \wedge f^{*}(\mathbb{T}_{F}) \longrightarrow \OO_{X}[n]$.
 By definition, $h$ gives us 
a \emph{homotopy} between $0$ and the composite morphism
$$\xymatrix{\mathbb{T}_{X} \wedge \mathbb{T}_{X} \ar[r] & 
f^{*}(\mathbb{T}_{F}) \wedge f^{*}(\mathbb{T}_{F}) \ar[r] & \OO_{X}[n].}$$
We let $\mathbb{T}_{f}$ be the relative tangent complex of $f$, so that we
 have an exact sequence of perfect complexes on $X$
$$\mathbb{T}_{f} \longrightarrow \mathbb{T}_{X} \longrightarrow f^{*}(\mathbb{T}_{F}).$$
The isotropic structure $h$ induces also a homotopy between $0$ and the composite morphism
$$\xymatrix{\mathbb{T}_{f} \otimes \mathbb{T}_{X} \ar[r] & \mathbb{T}_{X} \wedge \mathbb{T}_{X} \ar[r] & 
f^{*}(\mathbb{T}_{F}) \wedge f^{*}(\mathbb{T}_{F}) \ar[r] & \OO_{X}[n].}$$
As the morphism
$\mathbb{T}_{f} \longrightarrow f^{*}(\mathbb{T}_{F})$ comes itself with 
a canonical homotopy to $0$, by composing these homotopies, we end up with a loop pointed at $0$ in the space \linebreak
$Map_{L_{qcoh}(X)}(\mathbb{T}_{f} \otimes \mathbb{T}_{X},O_{X}[n])$. This loop defines
an element in 
$$\pi_{1}(Map_{L_{qcoh}(X)}(\mathbb{T}_{f} \otimes \mathbb{T}_{X},\OO_{X}[n]),0) \simeq
[\mathbb{T}_{f} \otimes \mathbb{T}_{X},\OO_{X}[n-1]].$$
By adjunction, we get a morphism of perfect complexes on $X$
$$\Theta_{h} : \mathbb{T}_{f} \longrightarrow \mathbb{L}_{X}[n-1],$$
depending on the isotropic structure $h$. 

\begin{df}\label{d11}
Let $f : X \longrightarrow F$ be a morphism of derived Artin stacks and $\omega$ an $n$-shifted symplectic form on $F$. An isotropic 
structure $h$ on $f$ is a \emph{Lagrangian structure on $f$ (relative to $\omega$)} if
the induced morphism
$$\Theta_{h} : \mathbb{T}_{f} \longrightarrow \mathbb{L}_{X}[n-1]$$
is a quasi-isomorphism of perfect complexes. 
\end{df}

The usefulness of Lagrangian structures is shown by the following existence theorem.

\begin{thm}\label{t2}
Let 
$$\xymatrix{
 & Y \ar[d]^-{g} \\
 X \ar[r]_-{f} & F,}$$
be a diagram of derived Artin stacks, $\omega \in \Symp(F,n)$ 
an $n$-shifted symplectic form on $F$, 
and $h$ (respectively, $k$) be
a Lagrangian structure on $f$ (respectively, on $g$). Then, the derived Artin stack 
$X\times_{F}^{h}Y$ is equipped with a canonical $(n-1)$-shifted symplectic structure 
called the \emph{residue} of $\omega$, and denoted by $R(\omega,h,k)$. 
\end{thm}

\textbf{Proof.} Let $Z:=X\times_{F}^{h}Y$. The two morphisms
$$p : \xymatrix{Z \ar[r] & X \ar[r] & F} \qquad q : \xymatrix{Z \ar[r] & Y \ar[r] & F}$$
come equipped with a natural homotopy $u : p \Rightarrow q$. This $u$ gives rise to 
a homotopy between the induced morphisms on the spaces of closed $2$-forms
$$u^{*} : p^{*} \Rightarrow q^{*} : \A^{2,cl}(X,n) \longrightarrow 
\A^{2,cl}(Z,n).$$
Moreover, $h$ (respectively, $k$) defines a path in the space $\A^{2,cl}(Z,n)$
$$h : 0 \leadsto p^{*}(\omega)$$  
(respectively, $$k : 0 \leadsto q^{*}(\omega).\, )$$
By concatenation of $h$, $u^{*}(\omega)$ and $k^{-1}$, we get 
a loop at $0$ in the space $\A^{2,cl}(Z,n)$, therefore a well
defined element
$$R(\omega,h,k) \in \pi_{1}(\A^{2,cl}(Z,n))\simeq \pi_{0}(\A^{2,cl}(Z,n-1)).$$
It remains to show that this closed $2$-form of degree $(n-1)$ is non-degenerate. This follows
from the definition of a Lagrangian structure. To see this, let 
$\pi:=p : Z \longrightarrow F$ be the natural map\footnote{Since we just want to prove that the closed $2$-form is non-degenerate, we might equivalently have chosen to run the argument for $\pi:= q$, instead of $\pi:=p$.}, and 
$$pr_{X} : Z \longrightarrow X \qquad pr_{Y} : Z \longrightarrow Y$$ the two projections.
By definition we have a commutative diagram
in $L_{qcoh}(Z)$ with exact rows
$$\xymatrix{\mathbb{T}_{Z} \ar[d]_-{\Theta_{R(\omega,h,k)}} \ar[r] & 
pr_{X}^{*}(\mathbb{T}_{X}) \oplus pr_{Y}^{*}(\mathbb{T}_{Y}) \ar[r] \ar[d]_-{\Theta_{h} \oplus
\Theta_{k}} & \pi^{*}(\mathbb{T}_{F}) \ar[d]^-{\Theta_{\omega}} \\
\mathbb{L}_{Z}[n-1] \ar[r] & pr_{X}^{*}(\mathbb{L}_{f})[n-1] \oplus pr_{Y}^{*}(\mathbb{L}_{g})[n-1] \ar[r]
& \pi^*(\mathbb{L}_{F}[n]).}$$
The morphism $\Theta_{\omega}$ is a quasi-isomorphism because $\omega$ is non-degenerate. 
The morphism $\Theta_{h} \oplus
\Theta_{k}$ is also a quasi-isomorphism because of the definition of Lagrangian structures. This
implies that $\Theta_{R(\omega,h,k)}$ is a quasi-isomorphism, and thus that 
$R(\omega,h,k)$ is an $(n-1)$-shifted symplectic structure. 
\hfil $\Box$ \\

An immediate corollary is the following statement.

\begin{cor}\label{ct2}
Let $X$ be a smooth Deligne-Mumford stack over $\mathrm{Spec}\, k$, $\omega$ a symplectic
form on $X$, and $L,L' \subset X$ two smooth closed Lagrangian substacks
(in the sense that $\omega$ vanishes on $L$ and $L'$, and both $L$ and $L'$ are of middle 
dimension). Then, the derived fiber product $L\times_{X}^{h}L'$
carries a canonical $(-1)$-shifted symplectic structure.
\end{cor}

\textbf{Proof.} It follows by the following simple observations. As $X$, $L$ and $L'$ are
smooth, the spaces of closed $2$-forms of degree $0$ on $X$, $L$ and $L'$ are (homotopically) discrete
(and equal to $H^{0}(X,\Omega^{2,cl}_{X/k})$,  $H^{0}(L,\Omega^{2,cl}_{L/k})$, 
$H^{0}(L',\Omega^{2,cl}_{L'/k})$). From this it follows that 
the spaces of isotropic and Lagrangian structures on the two inclusions
$$L \hookrightarrow X \qquad L' \hookrightarrow X$$
are either empty or (equivalent to) a point. As $L$ and $L'$ are Lagrangian substacks it is
easy to see that these spaces are non-empty and hence both equivalent to a point. In particular, there are
unique Lagrangian structures on the above two inclusion morphisms, and thus Theorem 
\ref{t2} implies that $L\times_{X}^{h}L'$ is endowed with a canonical 
$(-1)$-shifted symplectic structure.
\hfill $\Box$ \\

A particular case of the above corollary is the existence of $(-1)$-shifted symplectic forms on 
derived critical loci of functions on smooth Deligne-Mumford stacks (see also 
\cite{vez} for a more direct local approach). 

\begin{cor}\label{ct2'}
Let $X$ be a smooth Deligne-Mumford stack over $\mathrm{Spec}\, k$, 
and $f \in \OO(X)$ a global function on $X$, with differential
$df : X \longrightarrow T^{*}X$. Then, the derived critical locus of $f$, defined
as the derived fiber product
$$\mathbb{R}Crit(f):=X \times_{df,T^{*}X,0}X,$$
of the zero section with the section $df$ inside the total cotangent stack $T^{*}X$, carries
a canonical $(-1)$-shifted symplectic structure.
\end{cor}

\textbf{Proof.} We simply observe that $T^{*}X$ carries a canonical symplectic structure
and that $X$ sits inside $T^{*}X$ as two Lagrangian substacks, either via the zero section
of via the section $df$. \hfill $\Box$ \\

\subsection{$2$-shifted symplectic structure on $\Parf$}

We now state our third existence theorem, giving a canonical 
$2$-shifted symplectic structure on the derived stack of perfect complexes $\Parf$. This
$2$-shifted symplectic form will be constructed using the Chern character of the universal
object on $\Parf$, with values in negative cyclic homology. We will use the 
construction in \cite{chern}, as it is perfectly suited to our context, but any
functorial enough construction of the Chern character could be used instead. \\

We recall from \cite{tv} the definition of the derived stack $\Parf$, which was denoted there by 
$\mathcal{M}_{\mathbf{1}}$. The functor $\Parf$ sends a differential non-positively graded algebra $A$ to the nerve of the category of perfect (i.e. homotopically finitely presentable, or equivalently, dualizable in the monoidal model category of $A$-dg-modules) $A$-dg-modules which are cofibrant in the projective model structure of all $A$-dg-modules. It is a locally geometric derived stack, that is
a union of open substacks which are derived Artin stacks of finite presentation over $\Spec\, k$. 
Everything we said about derived Artin stacks locally of finite presentation also make
sense for $\Parf$, in particular we can speak about $p$-forms, closed $p$-forms and
symplectic structures on $\Parf$, even though $\Parf$ is not strictly speaking 
a derived Artin stack. 

On the derived stack $\Parf$ we have the universal perfect complex $\mathcal{E} \in L_{parf}(\Parf)$. 
The endomorphisms of this perfect complexes define a perfect dg-algebra 
over $\Parf$, denotes by
$$\mathcal{A}:=\rch (\mathcal{E},\mathcal{E}) \simeq \mathcal{E}^{\vee} \otimes \mathcal{E}.$$
The derived loop stack of $\Parf$
$$\mathcal{L}\Parf:=\mathbf{Map}(S^{1},\Parf) \longrightarrow \Parf$$
can be identified with the derived group stack over $\Parf$ of invertible elements in $\mathcal{A}$
$$\mathcal{L}\Parf \simeq \mathcal{A}^{*}:=GL_{1}(\mathcal{A}),$$
where $GL_{1}(\mathcal{A})$ is the group stack of auto-equivalences of
$\mathcal{A}$ considered as an $\mathcal{A}$-module. As $\mathcal{A}$ is 
a perfect dg-algebra over $\Parf$, we see that $\mathcal{A}^{*}$ is
a derived Artin group stack over $\Parf$, and that the corresponding 
quasi-coherent dg-Lie algebra over $\Parf$ is $\mathcal{A}$ itself, endowed
with its natural bracket structure given by the commutator. 
This implies that we have a natural quasi-isomorphism of
perfect complexes over $\Parf$
$$\mathbb{T}_{\Parf} \simeq \mathcal{A}[1].$$
This identification can be used to define a $2$-form of degree $2$ on $\Parf$
$$\mathbb{T}_{\Parf} \wedge \mathbb{T}_{\Parf}\simeq Sym^{2}(\mathcal{A})[2] 
\longrightarrow \OO_{\Parf}[2],$$
which, by definition, is (the shift by $2$ of) the composition
$$\xymatrix{\mathcal{A} \otimes \mathcal{A} \ar[r]^-{mult} & 
\mathcal{A} \ar[r]^-{Tr} & \OO_{\Parf}},$$
of the multiplication and the trace (or evaluation) morphism
$$\mathcal{A} \simeq E^{\vee}\otimes E \longrightarrow \OO_{\Parf}.$$
This $2$-form is clearly non-degenerate.
We will now see that this $2$-form of degree $2$ on $\Parf$ is the underlying
$2$-form of a canonical $2$-shifted symplectic structure on $\Parf$. \\

We consider two derived stacks
$$|NC| : \dAff_{k}^{op} \longrightarrow \mathbb{S} \qquad
HC^{-} : \dAff_{k}^{op} \longrightarrow \mathbb{S},$$
defined as follows. We have the derived stack 
in mixed complexes $\mathbf{DR}$ (we forget the extra grading here), 
on which we can apply the construction $NC$ to get a derived stack
in complexes of $k$-modules. The derived stack 
$|NC|$ is obtained by applying the $\s$-functor $E \mapsto |E|$ (i.e.
the Dold-Kan construction) to turn $NC$ into a derived stack 
of spaces. In other words, 
$|NC|$ sends a commutative $k$-dg-algebra $A$ to 
the simplicial set obtained from the complex $NC(\mathbf{DR}(A))$. 

The derived stack $HC^{-}$ sends an affine derived scheme $X$ to 
$\mathcal{O}(\mathcal{L}X)^{hS^{1}}$, the space of 
$S^{1}$-homotopy fixed functions on the loop space of $X$. It is easy to
describe this derived stack algebraically using simplicial commutative
algebras. To a commutative $k$-dg-algebra $A$ we form 
$A'$ the corresponding commutative simplicial $k$-algebra (see \cite{tove}), and
consider $S^{1}\otimes_{k}^{\mathbb{L}}A'$, which is another simplicial
commutative algebra on which the simplicial group $S^{1}=B\mathbb{Z}$ acts. 
The space of homotopy fixed point of this action is a model for
$HC^{-}(\Spec\, A)$
$$HC^{-}(\Spec\, A) \simeq (S^{1}\otimes_{k}^{\mathbb{L}}A')^{hS^{1}}.$$

The main theorem of \cite{tove} states that these two derived stacks
$|NC|$ and $HC^{-}$ are naturally equivalent (and that there is
moreover a unique equivalence respecting the multiplicative structures).
Furthermore, the Chern character construction of \cite[\S 4.2]{chern} produces a morphism of derived stacks
$Ch : \Parf \longrightarrow HC^{-}$, 
and thus, by the mentioned equivalence, a morphism
$$Ch : \Parf \longrightarrow |NC|.$$
This is the Chern character of the universal perfect complex, and defines a natural element in 
$$Ch(\mathcal{E}) \in \pi_{0}(Map(\Parf,|NC|)) \simeq H^{0}(NC(\Parf)).$$
We can project this element by the projection on the weight $2$ piece
$NC \longrightarrow NC^{w}(2)$ to get 
$$Ch(\mathcal{E})_{2} \in H^{0}(NC^{w}(\Parf)(2)) \simeq \pi_{0}(\A^{2,cl}(\Parf,2)),$$
which is a closed $2$-form of degree $2$. 

\begin{thm}\label{t3}
The closed $2$-form 
$Ch(\mathcal{E})_{2}$ defined above is a $2$-shifted symplectic structure on $\Parf$. 
\end{thm}

\textbf{Proof.} The proof consists in identifying the underlying 
$2$-form (of degree 2) of $Ch(\mathcal{E})_{2}$, and show that it coincides  (up to a factor $\frac{1}{2}$) 
with 
the $2$-form described earlier in this section
$$\xymatrix{\mathcal{A} \otimes \mathcal{A} \ar[r]^-{mult} & 
\mathcal{A} \ar[r]^-{Tr} & \OO_{\Parf}}.$$
This identification can be seen as follows. For a perfect complex $E$ on a derived Artin stack $X$, 
its Chern character $Ch(E)$ has an image in Hodge cohomology
$$Ch(E)=\sum Ch_{p}(E) \in \oplus_{p} H^{p}(X,\wedge^{p}\mathbb{L}_{X/k}),$$
induced by the natural morphism from negative cyclic homology to Hochschild
homology (i.e. the projection $NC \rightarrow |\mathbf{DR}|$). 
These Hodge cohomology classes can be described using the Atiyah class of $E$
$$a_{E} : E \longrightarrow E\otimes_{\OO_{X}}\mathbb{L}_{X/k}[1].$$
We can compose this class with itself to get
$$a_{E}^{i} : E \longrightarrow E\otimes_{\OO_{X}} (\wedge^{i}\mathbb{L}_{X/k})[i],$$
which we write as
$$a_{E}^{i} : \OO_{X} \longrightarrow E^{\vee}\otimes_{\OO_{X}} E \otimes_{\OO_{X}} 
(\wedge^{i}\mathbb{L}_{X/k})[i].$$
Composing with the trace $E^{\vee}\otimes_{\OO_{X}} E \longrightarrow \OO_{X}$, we obtain
classes in $H^{i}(X,\wedge^{i}\mathbb{L}_{X/k})$. We have
$$Ch_{i}(E)=\frac{Tr(a_{E}^{i})}{i !} \in H^{i}(X,\wedge^{i}\mathbb{L}_{X/k}).$$
We come back to our specific situation where $X=\Parf$ and $E=\mathcal{E}$ is the
universal perfect complex. The shifted cotangent complex 
$\mathbb{L}_{X/k}[1]$ is naturally equivalent to $\mathcal{E}\otimes_{\OO_{X}}\mathcal{E}^{\vee}$, and the
Atiyah class 
$$a_{\mathcal{E}} : \mathcal{E} \longrightarrow \mathcal{E}\otimes_{\OO_{X}}\mathcal{E}\otimes_{\OO_{X}}\mathcal{E}^{\vee}$$
is simply the adjoint of the multiplication morphism
$$\mathcal{E} \otimes \mathcal{A} \longrightarrow \mathcal{E},$$
where $\mathcal{A}=\rch(\mathcal{E},\mathcal{E})$ is the endomorphism dg-algebra
of $\mathcal{E}$. From this we get the required formula 
$$Ch_{2}(\mathcal{E})=\frac{1}{2}\cdot Tr(mult) \in H^{2}(\Parf,\wedge^{2}\mathbb{L}_{\Parf}).$$
\hfill $\Box$ \\

As a corollary of the last theorem and of theorem \ref{t1}, we have the 
following statement, which is an extension of the corollary \ref{ct1}
from the case of vector bundles to the case of perfect complexes. 

\begin{cor}\label{ct3}
Let $Y$ be a smooth and proper
Deligne-Mumford stack with connected geometric fibers of relative dimension $d$.
\begin{enumerate}
\item The choice of a fundament class $[Y] \in H^{2d}_{DR}(Y,\OO)$
determines a canonical $2(1-d)$-shifted symplectic form on the derived stack 
of perfect complexes with flat connexions on $Y$ 
$$\Parf_{DR}(Y):=\mathbf{Map}(Y_{DR},\Parf).$$
\item The choice of a fundament class $[Y] \in H^{2d}_{Dol}(Y,\OO)$
determines a canonical $2(1-d)$-shifted symplectic form on the derived stack 
of perfect complexes with Higgs fields
$$\Parf_{Dol}(Y):=\mathbf{Map}(Y_{Dol},\Parf).$$
\item The choice of a trivialization (when it exists)
$\omega_{Y/k} \simeq \OO_{Y}$,  
determines a canonical $2-d$-shifted symplectic form on the derived stack 
of perfect complexes on $Y$ 
$$\Parf(Y):=\mathbf{Map}(Y,\Parf).$$
\item If $M$ is a compact, orientable topological manifold of dimension $d$, then a choice
of a fundamental class $[M]\in H_{d}(M,k)$ determines a canonical 
$(2-d)$-shifted symplectic form on the derived stack of 
perfect complexes on $M$
$$\Parf(M):=\mathbf{Map}(M,\Parf).$$
\end{enumerate}
\end{cor}

\section{Examples and applications}

We present in this last section some examples and consequences of our results.

\subsection{$0$-shifted symplectic structures on moduli of sheaves on curves and surfaces}

We start here by explaining how theorems \ref{t1} and \ref{t3} can be used in order to 
recover the existence of well known symplectic forms on certain moduli spaces 
of local systems on curves, and vector bundles (or more generally of perfect complexes)
on K3 and abelian surfaces. \\

\noindent \textbf{Local systems on curves.} To start with, assume that $G$ is a simple algebraic group over some field $k$, and 
$C$ be a (geometrically connected) smooth and proper curve over $k$. We have the
derived moduli stacks $\mathbb{R}\textrm{Loc}_{DR}(C,G)$, $\mathbb{R}\textrm{Loc}_{Dol}(C,G)$
and $\mathbb{R}\textrm{Loc}_{B}(C,G):=\mathbb{R}\textrm{Loc}(C^{top},G)$, of local systems of principal $G$-bundles on $C$, 
Higgs $G$-bundles on $C$, and flat $G$-bundles on the underlying topological space $C^{top}$
of $C$. According to our corollary \ref{ct1}, a choice of an orientation of $C$ determines
$0$-shifted symplectic structures on these spaces. These spaces contain
smooth Deligne-Mumford substacks consisting of simple objects
$$\mathrm{Loc}_{DR}(C,G)^{s} \subset \mathbb{R}\textrm{Loc}_{DR}(C,G) \qquad
\textrm{Loc}_{Dol}(C,G)^{s} \subset \mathbb{R}\textrm{Loc}_{Dol}(C,G)$$ 
$$\textrm{Loc}_{B}(C,G)^{s} \subset \mathbb{R}\textrm{Loc}_{B}(C,G).$$
These substacks of simple objects are moreover \'etale gerbes over 
smooth algebraic varieties, bounded by the center of $G$ (so they are
algebraic varieties as soon as this center is trivial). Therefore, 
the restriction of these $0$-shifted symplectic forms on these substacks define
symplectic forms on the corresponding 
smooth algebraic varieties. We recover this way well known symplectic structures
on the coarse moduli space of simple flat $G$-bundles, simple Higgs $G$-bundles, and
simple flat $G$-bundles on $C^{top}$ (see e.g. \cite{go,je,iis}). 

It is interesting to note here that 
our results imply that these symplectic forms existing on the smooth locus
of simple objects have canonical extension to the whole derived moduli stacks. Another
interesting remark is the case of $C=\mathbb{P}^{1}$ is a projective line.
The corresponding coarse moduli of simple objects is just a point, but the
derived stacks are non-trivial. For instance, in the de Rham setting, 
we have
$$\mathbb{R}\textrm{Loc}_{DR}(C,G) \simeq [(\Spec\, Sym(\mathfrak{g}^{\vee}[1]))/G].$$
This derived stack has already been considered in \cite{laf}, and according to our
results carries a canonical $0$-shifted symplectic form. The tangent complex at the
unique closed point is
$$\mathbb{T} \simeq \mathfrak{g}[1] \oplus \mathfrak{g}[-1],$$
and the $0$-shifted symplectic form there is the canonical 
identification
$$\mathfrak{g}[1] \oplus \mathfrak{g}[-1] \simeq \mathfrak{g}^{\vee}[1] \oplus \mathfrak{g}^{\vee}[-1]$$
induced by the isomorphism $\mathfrak{g} \simeq \mathfrak{g}^{\vee}$ given by 
the data of the (essentially unique) symmetric, non-degenerate, bilinear $G$-invariant
form on $\mathfrak{g}$. \\

\noindent \textbf{Perfect complexes on CY surfaces.} 
As a second example we let $S$ be a Calabi-Yau surface (either a K3 or an abelian surface)
over a field $k$, equipped with a trivialization $\omega_{S}\simeq \OO_{S}$. 
We have $\Parf(S)$, the derived moduli stack of perfect complexes
on $S$. According to our corollary \ref{ct3}, this derived stack is equipped with 
a canonical $0$-shifted symplectic form. We let $\Parf(S)^{s} \subset \Parf(S)$
be the derived open substack of simple objects, that is perfect complexes
with no negative self extensions and only scalar multiplication as endomorphisms
(see \cite{tv}). We denote by $\mathcal{M}_{S}^{s}$ its truncation
$$\mathcal{M}_{S}^{s}:=h^{0}(\Parf(S)^{s}).$$
The (underived) stack $\mathcal{M}_{S}^{s}$ is a $\mathbb{G}_{m}$-gerbe over
an algebraic space $M_{S}^{s}$, locally of finite presentation over $k$.
It is proven in \cite{in} that $M_{S}^{s}$ is smooth and comes equipped with 
a natural symplectic structure. The existence of a non-degenerate $2$-form 
is easy, but there are relatively heavy computations in order to prove
that this $2$-form is closed. We will explain how this symplectic
structure can be deduced from the $0$-shifted symplectic structure 
on the whole derived stack $\Parf(S)$. First of all, 
the $0$-shifted symplectic form restricts to a $0$-shifted symplectic form
on the open $\Parf(S)^{s}$. We consider the 
determinant morphism of \cite[3.1]{stv}
$$det : \Parf(S)^{s} \longrightarrow \mathbb{R}Pic(S),$$
where $\mathbb{R}Pic(S)$ is the derived Picard stack of $S$, defined 
to be $\textbf{Map}_{\dSt_{k}}(S,B\mathbb{G}_{m})$. As explained in 
\cite[4.2]{stv} there is a natural projection
$\mathbb{R}Pic(S) \longrightarrow \Spec\, k[e_{1}]$,
with $e_{1}$ of degree $-1$. The choice of a point $s\in S(k)$ (assume there is one
for simplicity) defines by pull-back along $s : \Spec\,k \longrightarrow S$, 
another projection
$\mathbb{R}Pic(S) \longrightarrow \mathbb{R}Pic(\Spec\, k)=B\mathbb{G}_{m}$.
These two projections, pre-composed with the determinant map, defines a morphism of derived stacks
$$\pi : \Parf(S)^{s} \longrightarrow 
\Spec\, k[e_{1}] \times B\mathbb{G}_{m}.$$
It is easy to see that this projection is smooth and representable by an algebraic space. Thus, 
its fiber at the natural base point is a smooth algebraic space $X$, equipped with a natural
morphism
$$j : X \longrightarrow \Parf(S)^{s}.$$
By definition, $X$ is naturally isomorphic to $M_{S}^{s}$, the coarse moduli of
the truncation $\mathcal{M}_{S}^{s}$ of $\Parf(S)^{s}$.  Finally, 
the $0$-shifted symplectic form on $\Parf(S)^{s}$ can be pulled back 
to $X$ by $j$, and defines a closed $2$-form on $X$. The tangent 
of $X$ at a point corresponding to a perfect complex $E$ is $Ext^{1}(E,E)$, and this
$2$-form is the natural pairing
$$\xymatrix{Ext^{1}(E,E) \times Ext^{1}(E,E) \ar[r]^-{\cap} & 
Ext^{2}(E,E) \ar[r]^-{tr} & Ext^{2}(\OO_{S},\OO_{S}) \simeq k},$$
and thus is a symplectic form on the smooth algebraic space $X$, which 
is the one constructed in \cite{in}.

Again, the interesting remark here is that this symplectic structure on 
$M_{S}^{s}$ is induced by a $0$-shifted symplectic structure on the \emph{whole
derived stack} $\Parf(S)$, which includes all perfect complexes,
and in particular non-simple ones, and even complexes with possibly non trivial 
negative self-extensions.

\begin{rmk}\label{r4}
As the morphism 
$$\pi :  \Parf(S)^{s} \longrightarrow \Spec\, k[e_{1}] \times B\mathbb{G}_{m}$$
is smooth with fibers $X$, it is locally (for the smooth topology) on $\Parf(S)^{s}$
equivalent to the projection
$$X \times \Spec\, k[e_{1}] \times B\mathbb{G}_{m} \longrightarrow \Spec\, k[e_{1}] \times B\mathbb{G}_{m}.$$
To be more precise, we have a commutative diagram of derived stacks, with
cartesian squares
$$\xymatrix{
\Parf(S)^{s} \ar[r]^-{\pi} & \Spec\, k[e_{1}] \times B\mathbb{G}_{m} \\
Y \ar[u]^-{p} \ar[r] & \Spec\, k[e_{1}] \ar[u] \\
X \ar[u]^-{i}  \ar[r] & \Spec\, k. \ar[u]}$$
Here, $Y$ and $X$ are derived algebraic spaces, respectively smooth over 
$\Spec\, k[e_{1}]$ and $\Spec\, k$. Moreover, $X\simeq h^{0}(Y)$ is the truncation of $Y$, 
and $i : X \rightarrow Y$ is the natural closed embedding. Locally for the \'etale 
topology, $Y$ is a direct product $U \times \Spec\, k[e_{1}]$, where
$U$ is an \'etale scheme over $X$. Indeed, any affine scheme $Z$ smooth over
$\Spec\, k[e_{1}]$, splits (uniquely , the space of splitting is
connected) as a product $Z_{0} \times \Spec\, k[e_{1}]$ where 
$Z_{0}$ is a smooth affine variety. Finally, the morphism $p$ is 
a $\mathbb{G}_{m}$-torsor and thus is smooth and surjective. This shows that 
locally for the smooth topology on $\Parf(S)^{s}$, the morphism $\pi$ is
a direct product. 
\end{rmk}

\begin{rmk}\label{r4'}
It is also possible to show that this local decomposition of the previous
remark \ref{r4} is compatible with
the symplectic structure. More precisely, locally $\Parf(S)^{s}$ looks like
$X \times \Spec\, k[e_{1}] \times B\mathbb{G}_{m}$. Now, the derived
stack $\Spec\, k[e_{1}] \times B\mathbb{G}_{m}$ possesses a natural 
$0$-shifted symplectic structure, for instance by using our corollary \ref{ct1}
and the fact that we have a natural identification
$$\Spec\, k[e_{1}] \times B\mathbb{G}_{m} \simeq \Map(S^{2},B\mathbb{G}_{m}).$$
The projection $\pi$ is a symplectic morphism, and the local decomposition
$$\Parf(S)^{s} \simeq X \times \Spec\, k[e_{1}] \times B\mathbb{G}_{m}$$
becomes a decomposition of symplectic derived stacks. 
\end{rmk}

\subsection{$(-1)$-shifted symplectic structures and symmetric obstruction theories}

We compare here our notion of $(-1)$-shifted symplectic structure with the notion of 
a symmetric obstruction theory of \cite[Def. 1.10]{befan}. Recall first from \cite[\S 1]{stv} that for any 
derived stack $F$ which is locally of finite presentation, its
truncation $h^{0}(F)$ comes equipped with a natural perfect obstruction theory. 
It is constructed by considering the inclusion $j : h^{0}(F) \longrightarrow F$, and
by noticing that the induced morphism of cotangent complexes
$$j^{*} : j^{*}(\mathbb{L}_{F/k}) \longrightarrow \mathbb{L}_{h^{0}(F)}$$
satisfies the property to be a perfect obstruction theory. In practice all 
obstruction theories arise this way.

Assume now that $F$ comes equipped with a $(-1)$-shifted symplectic structure $\omega$. We write
the underlying $2$-form of degree $-1$ as a morphism of perfect complexes
$$\omega : \mathbb{T}_{F} \wedge \mathbb{T}_{F} \longrightarrow \OO_{F}[-1].$$
We use the fact that $\omega$ is non-degenerate and the equivalence
$$\Theta_{\omega} : \mathbb{T}_{F} \simeq \mathbb{L}_{F}[-1],$$
to get another morphism
$$Sym^{2}(\mathbb{L}_{F})[-2] \simeq (\mathbb{L}_{F}[-1]) \wedge 
(\mathbb{L}_{F}[-1]) \simeq  \mathbb{T}_{F} \wedge \mathbb{T}_{F} \longrightarrow \OO_{F}[-1],$$
which we rewrite as
$$Sym^{2}(\mathbb{L}_{F}) \longrightarrow \OO_{F}[1].$$
This pairing stays non-degenerate, and thus defines a equivalence
$$\mathbb{L}_{F} \simeq \mathbb{T}_{F}[1].$$
When restricted to the truncation $h^{0}(F)\hookrightarrow F$, we find that the perfect obstruction
theory $E:= j^{*}(\mathbb{L}_{F/k})$, comes equipped with a natural equivalence
$$E \simeq E^{\vee}[1],$$
which is symmetric (i.e. comes from a morphism $Sym^{2}(E) \rightarrow \OO[1]$). It is, by 
definition \cite[Def. 1.10]{befan}, a symmetric obstruction theory on $h^{0}(F)$. 

\begin{rmk}\label{r5}
As explained in \cite{stv}, the datum of the truncation $h^{0}(F)$ together with the
obstruction theory $j^{*}(\mathbb{L}_{F})$ is strictly weaker than the datum of  
the derived stack $F$. In the same way, the passage from a 
$(-1)$-shifted symplectic form on $F$ to a symmetric obstruction theory 
on $h^{0}(F)$ looses important informations. The most important one is that 
the corresponding symmetric obstruction theory only depends on the
underlying $2$-form of $\omega$, and thus does not see that $\omega$
comes with the important further closedness datum.\\
An interesting related question is whether a symmetric obstruction theory that is induced by a $(-1)$-shifted symplectic form, is \'etale locally isomorphic to the canonical one given by lagrangian intersections on a smooth scheme, or even to that existing on 
the derived zero locus of a \emph{closed} $1$-form on a smooth scheme - instead of just an \emph{almost-closed} $1$-form, as in the case of a general symmetric obstruction theory, see \cite[\S 3.4]{behrend}. Some very interesting formal and local (for the analytic topology) results in this direction have been proven in \cite{joyce} slightly after the appearance of the first preprint version of this paper, and efficiently applied by the same authors to show the existence of local potentials for Donaldson-Thomas theory on Calabi-Yau $3$-folds. 
\end{rmk}

Three main sources of examples of $(-1)$-shifted symplectic structures, and thus of
symmetric obstruction theories, are the following. \\

\noindent \textbf{Sheaves on CY 3-folds.} Let $X$ be a smooth and proper CY manifold of
dimension $3$, together with a trivialization $\omega_{X/k}\simeq \OO_{X}$. 
It is a $\OO$-compact object  endowed with an $\OO$-orientation of dimension $3$, so
our theorem \ref{t1} can be applied. We find in particular that the
derived stack $\Parf(X)$ of perfect complexes on $X$ is canonically endowed 
with a $(-1)$-shifted symplectic structure. This defines a symmetric obstruction
theory on the truncation $h^{0}(\Parf(X))$. The $(-1)$-shifted symplectic structure also
induces a $(-1)$-shifted symplectic structure on the derived stack of perfect complexes
with fixed determinant. For this, we use the determinant map of \cite{stv}
$$det : \Parf(X) \longrightarrow \mathbb{R}Pic(X),$$
and consider the fiber at a given global point $L \in Pic(X)$, corresponding to a line bundle on $X$
$$\Parf(X)_{L}:=det^{-1}(\{L\}).$$
The $(-1)$-shifted symplectic form on $\Parf(X)$ can be pulled-back to a closed $2$-form
on $\Parf(X)_{L}$ using the natural morphism $\Parf(X)_{L} \longrightarrow \Parf(X)$. 
It is easy to see that this closed $2$-form stays non-degenerate, and thus
defines a $(-1)$-shifted symplectic structure on $\Parf(X)_{L}$. Finally, restricting 
to simple objects, we get a quasi-smooth derived stack $\Parf(X)_{L}^{s}$, endowed
with a $(-1)$-shifted symplectic structure. On the truncation $h^{0}(\Parf(X)^{s}_{L})$ we thus find
a perfect symmetric obstruction theory of amplitude $[-1,0]$. 

We could also consider similar examples, like $\mathbb{R}\textrm{Loc}(M,G)$ - 
for $M$ an oriented compact 3-dimension topological manifold, and $G$ a reductive group
scheme over $k$ - or $\Parf(M):=\Map(M,\Parf)$. \\

\noindent \textbf{Maps from elliptic curves to a symplectic target.} The second source of examples
of $(-1)$-shifted symplectic structures is by considering maps from an elliptic curve towards
a symplectic smooth target. Let $E$ be a fixed elliptic curve endowed with 
a trivialization $\omega_{E/k}\simeq \OO_{E}$, and $(X,\omega)$ a smooth 
symplectic variety over $k$ (for instance $X$ can be $T^{*}Z$ for some
smooth $Z$, with its canonical symplectic form. This example is of
fundamental interests in elliptic cohomology, see \cite{cg}). We let 
$F:=\Map(E,X)$, which by our theorem \ref{t1} is endowed with 
a canonical $(-1)$-shifted symplectic structure. The derived stack $F$ is a derived algebraic 
space, and in fact a quasi-projective derived scheme if $X$ is itself quasi-projective. 
It is moreover quasi-smooth, as its tangent at a given point $f : E \rightarrow X$ is 
the Zariski cohomology complex $C^{*}(E,f^{*}(T_{X}))$. We thus have 
a quasi-smooth derived algebraic space (or even scheme), endowed with 
a $(-1)$-shifted symplectic structure. It gives rise to a symmetric perfect obstruction
theory of amplitude $[-1,0]$ on the truncation $h^{0}(F)$. Once again,
the existence of the $(-1)$-shifted symplectic form is a stronger statement than the existence
of a symmetric obstruction theory. This $(-1)$-shifted symplectic structure induces, by passing to the derived formal completion at a given point, the degree $-1$ symplectic structure recently considered by Costello in \cite{co}. As explained in \cite{co}
it can be used to construct a quantization of the moduli $F$, and 
for this the datum of the symmetric obstruction theory is not enough. \\

\noindent \textbf{Lagrangian intersections.} Let $(X,\omega)$ be a smooth 
symplectic scheme over $k$, with two smooth Lagrangian subschemes 
$L$ and $L'$. Then, the two closed immersions $L,L' \subset X$ are
endowed with a unique Lagrangian structure in the sense of our definition \ref{d11}, and thus
our theorem \ref{t2} implies that the derived intersection
$L\times^{h}_{X}L'$ carries a natural $(-1)$-shifted symplectic structure. Again this
defines a symmetric perfect obstruction theory of amplitude $[-1,0]$ on the
truncation, that is on the usual schematic fiber product $L\times_{X}L'$. The
data of the $(-1)$-shifted symplectic structure is again stronger than the data of the
corresponding symmetric obstruction theory, as for instance it can be used to 
\emph{quantize} the intersection $L\times_{X}L'$. 
We will come back to the quantization construction in a future work.

\end{document}